\documentclass[12pt]{amsart}
\usepackage{amsmath,amsthm,amsfonts,amssymb}
\begin{document} 
\newcommand{\A}{{\mathbb A}}
\newcommand{\B}{{\mathbb B}}
\newcommand{\C}{{\mathbb C}}
\newcommand{\N}{{\mathbb N}}
\newcommand{\Q}{{\mathbb Q}}
\newcommand{\Z}{{\mathbb Z}}
\renewcommand{\P}{{\mathbb P}}
\newcommand{\R}{{\mathbb R}}
\newcommand{\rc}{\subset}
\newcommand{\rank}{\mathop{rank}}
\newcommand{\tensor}{\mathop{\otimes}}
\newcommand{\trace}{\mathop{tr}}
\newcommand{\codim}{\mathop{codim}}
\newcommand{\dimc}{\mathop{dim}_{\C}}
\newcommand{\Lie}{\mathop{Lie}}
\newcommand{\Sing}{\mathop{Sing}}
\newcommand{\Spec}{\mathop{Spec}}
\newcommand{\Auto}{\mathop{{\rm Aut}_{\mathcal O}}}
\newcommand{\Aut}{\mathop{{\rm Aut}}}
\newcommand{\alg}[1]{{\mathbf #1}}
\newcommand{\coeff}{\mathop{{\text Coeff}}}
\newtheorem{lemma}{Lemma}[section]
\newtheorem{definition}[lemma]{Definition}
\newtheorem*{claim}{Claim}
\newtheorem{corollary}[lemma]{Corollary}
\newtheorem{Conjecture}[lemma]{Conjecture}
\newtheorem{assumptions}[lemma]{Assumptions}
\newtheorem*{SpecAss}{Special Assumptions}
\newtheorem{example}[lemma]{Example}
\newtheorem*{remark}{Remark}
\newtheorem*{observation}{Observation}
\newtheorem*{fact}{Fact}
\newtheorem*{remarks}{Remarks}
\newtheorem*{question}{Question}

\newtheorem{proposition}[lemma]{Proposition}
\newtheorem{theorem}[lemma]{Theorem}
\numberwithin{equation}{section}
\def\labelenumi{\rm(\roman{enumi})}
\newcommand{\XXX}{\message{XXX}.{\bf XXX !}}
\title{Dense entire curves in Rationally Connected manifolds
}
\author {Fr\'ed\'eric Campana \& J\"org Winkelmann}

\begin{abstract} We show the existence of metrically dense entire curves (of growth order $0)$ in rationally connected complex projective manifolds, confirming for this case
  a conjecture formulated in \cite{Ca04}, according to which such entire curves on projective manifolds exist if and only if these are `special' in the sense defined in loc.cit. 
  For unirational manifolds, the statement above is an easy consequence of function theoretic methods. Our proof rests on the `comb smoothing' technique of Koll\'ar-Miyaoka-Mori, and may be seen as a substitute of the power series expansion of entire functions with values in $\C^n$, in the absence of global coordinates.

   We next show the existence of dense entire curves, avoiding the singular locus, in certain log-terminal normal rational surfaces. This implies via results of Grassi  and Oguiso the existence of dense entire curves into any Calabi-Yau threefold fibered in Abelian surfaces or elliptic curves.

  We then show that a dense entire curve may be chosen on any rationally connected manifold in such a way that it does not lift
  to any of its ramified covers, answering in this case a question of \cite{CZ} about the Nevanlinna analog of the `weak Hilbert property' of arithmetic geometry. We consider briefly the other test case of the conjecture, namely manifolds with $c_1=0$, and most `special' surfaces.
  
\end{abstract}

\address{%
Fr\'ed\'eric Campana \\
D\'epartement de Math\'ematiques \\
Universit\'e Nancy 1 \\
Vandoeuvre-l\`es-Nancy\\
54500 \\
France
}

\email{frederic.campana@univ-lorraine.fr
}

\address{%
J\"org Winkelmann \\
IB 3-111\\
Lehrstuhl Analysis II \\
Fakult\"at f\"ur Mathematik \\
Ruhr-Universit\"at Bochum\\
44780 Bochum \\
Germany
}
\email{joerg.winkelmann@rub.de\newline
    ORCID: 0000-0002-1781-5842
}

\maketitle

\tableofcontents

   \section{Introduction}

   Our main result is the existence of dense entire curves
   (i.e.~a holomorphic map from $\C$ with metrically dense image) in rationally
connected manifolds.

More precisely, we prove that, given a rationally connected complex
projective manifold $X$ with a countable subset $M$, a
hypersurface $D$ and a non-trivially ramified cover $Y\to X$ we
can find a dense entire curve $c:\C\to X$ such that
\begin{itemize}
\item
  The countable set $M$ is contained in the image $c(\C)$,
\item
  $c(\C)$ meets $D$ tranversally somewhere.
\item
  $c$ can not be lifted to an entire curve in $Y$.
\item
  The order $\rho_c$ (in the sense of Nevanlinna theory)
  equals zero.
\end{itemize}

The non-liftability confirms (for rationally connected manifolds)
a conjecture of Corvaja and Zannier.

We also prove some results in these directions for normal projective
surfaces and for compact K\"ahler manifolds with $c_1=0$.

The existence of a dense entire curve for any rationally connected
manifolds confirms, for this case,
a conjecture, which motivated the present text, that for compact K\"ahler manifolds, the existence
of a dense entire curve is equivalent to being `special' (definition recalled below). Special manifolds are in a precise sense `opposite'
to manifolds of general type. While, following a conjecture of S. Lang, the non-existence of Zariski dense entire curves on manifolds of general 
type has been investigated since decades, our result  seems to be the first one in the opposite direction, beyond the classical case of unirational or abelian
manifolds.

Although such dense curves exist on any unirational manifold for simple function-theoretic interpolation reasons, it is presently unknown whether rationally connected are all unirational. In fact the contrary is expected to be true quite generally. Moreover our proof, which is mainly algebraic, based on the `comb smoothing'
technique of \cite{KMM}, provides a (to our knowledge) new approach to these interpolation properties, and plays a r\^ole similar to the power series expansion of entire functions in this broader context. We don't know however if any entire curve in a rationally connected manifold can be approximated by algebraic maps, as in our construction.

 The conjectural links formulated by S. Lang between arithmetic and hyperbolicity properties lead to conjecture that the existence of a Zariski dense entire curve should occur on $X$ defined over a number field $k$ if and only if $X(k')$ is Zariski dense for some finite extension $k'$ of $k$. As this arithmetic property is unknown for $X$ rationally connected, our result give some support to this arithmetic conjecture.

We show finally that  the entire curves we construct
in any rationally connected manifolds may be chosen so as to have growth order
$\rho_f=0$.

In particular, any rationally connected complex projective manifold $X$
admits dense entire curves $h:\C\to X$ with $\rho_h=0$. In the opposite direction
we showed in \cite{CW"} that any complex projective manifold
$X$ of dimension $n$ admitting a non-degenerate holomorphic
map $F:\C^n\to X$ of order $\rho_F<2$ is
rationally connected.

For normal projective surfaces $X$ with only quotient singularities we show:
If there is an effective non-zero $\Q$-divisor $\Delta$ such that $K+\Delta$
is $\Q$-trivial and such that the pair $(X,\Delta)$ is log terminal,
then there exists a dense entire curve avoiding the singular locus. Here the $2$-dimensional MMP plays a crucial r\^ole in the proof.

\section{Brief review of special manifolds}

Let $X$ be  an $n$-dimensional smooth connected compact complex manifold (either projective or compact K\"ahler). Its Kodaira dimension is denoted with $\kappa(X)$.

Recall (see \cite{Ca04}) that a compact K\"ahler manifold $X$ is said to be `special' if no rank one subsheaf $L\subset \Omega^p_X$ has top Kodaira dimension $p$, this for any $p>0$. These manifolds are higher dimensional generalisations of rational and elliptic curves. In particular, manifolds $X$ which are either rationally connected (see next section for a brief reminder), or with $\kappa(X)=0$ are special\footnote{However, among $n$-dimensional manifolds with $\kappa=k$,  there are both special and non-special examples, for any $k\neq 0,n$.}.  In fact, the decomposition $c= (j\circ r)$ of the `core map' (See \cite{Ca04} for details) shows that the `building blocks' of special manifolds should be  (smooth) `orbifold pairs' either `rationally connected'\footnote{in the sense of having $\kappa^+=-\infty$. When there is no orbifold structure, a classical conjecture claims that this is equivalent to rational connectedness.}, or with $\kappa=0$. Moreover (\cite{Ca04}, Corollary 8.11) $X$ is `special' if it is $\C^n$-dominable (i.e. if there exists a meromorphic map from $\C^n$ to $X$ regular and submersive at some point). The converse is not expected to hold in general (but partially weaker versions  are conjectured below).

We conjecture that the special manifolds are also characterised by the two (conjecturally equivalent properties) that their Kobayashi pseudometric $d_X$ vanishes identically, and that they admit a holomorphic map $h:\C\to X$ with metrically dense image. This conjecture is motivated by the above-mentioned decomposition $c=(j\circ r)^n$ of the core map, which essentially reduces these conjectures to the two cases of orbifold pairs either rationally connected in the above sense, or with $\kappa=0$.

For rationally connected manifolds, $d_X$ vanishes obviously. We prove the existence of a dense $h(\C)$ in Theorem \ref{rc-main}, using the comb-smoothing technique of \cite{KMM}. We also deduce a Nevanlinna analogue of the Hilbert Irreducibility Property introduced in \cite{CZ}. Similar results hold more classically when $X$ is a complex torus. 

We give some brief remarks on the much more challenging class of manifolds with zero first Chern class in section \ref{c1=0} below.

$\bullet$ This text was completed as the first named author (F.C) was staying at the Mathematical Sciences Research Institute in Berkeley, California, in march 2019, supported by the National Science Foundation under Grant No. 1440140. He thanks J. Mc Kernan and S. Kebekus for informations about rational curves on normal surfaces, used in Proposition \ref{rsurf} below. We also thank F. Forstneric for informing us about the text of A. Gournay, and J. Noguchi for asking us about the possible growth orders of the entire curves we construct, which led to the addition of the final section.

\section{Review of Rationally connected manifolds}

Recall that a projective manifold $X$ is said to be:

1. {\em rational} if it is bimeromorphic to $\P^n$, 

2. {\em unirational} if dominated by $\P^n$ (which means that there exists a non-degenerate
(i.e.: dominant) meromorphic map $\Phi:\P^n\dasharrow X)$ , and:

3. {\em rationally connected} if any two generic points of $X$ are connected by a rational curve in $X$.

4. {\em Fano} if its anticanonical bundle is ample.

These properties (except for 4) are birational.

One has the implications: 

\[
\begin{array}{ccccc}
  &&&&  \text{Fano}\\
  &&&&  \Downarrow \\
  \text{rational} & \Longrightarrow
  & \text{unirational}  & \Longrightarrow
  & \text{rationally connected}
  \\
\end{array}
\]

\smallskip

For curves (resp. surfaces), these properties (resp. except being Fano)
are equivalent. But, starting
in dimension $3$, many unirational manifolds are known to be non-rational. On the other hand, although it is expected that most rationally connected threefolds are not unirational, there is pre\-sent\-ly no known invariant to distinguish them.

Examples of rationally connected threefolds not known to be unirational, and possibly non-unirational, are `general' quartics in $\P^4$, double covers of $\P^3$ branched over a smooth sextic, and standard conic bundles over $\P^2$ with smooth discriminant of large degree.

\section{Conjectures on Special manifolds}

\begin{definition}(\cite{Ca04}, Definition 2.24 and Theorem 2.27) A compact K\"ahler manifold $X$ is said to be `special' if $\kappa(X,L)<p$ for any rank-one coherent sheaf $L\subset \Omega^p_X$, and any $p>0$.
\end{definition}

We refer to \cite{Ca04} for details on this class of manifolds, and some of the reasons to introduce them, the main one being the existence of the `core map' which splits any $X$ into its opposite parts: special vs general type. We shall here simply mention the two basic examples of `special' manifolds:
Rationally connected manifolds and manifolds with $\kappa=0$.

Below we summarize what we expect to be true for special manifolds. Many more variants can be formulated, we state only the simplest ones.
\begin{Conjecture}\label{conjspec}(\cite{Ca04},
  \S9.2, 9.8, 9.5, 9.20.)
  \begin{enumerate}
   \item
  A compact connected K\"ahler manifold $X$ is `special' if and only if any one of the following properties is satisfied:

1.1. The Kobayashi pseudo-metric of $X$ vanishes identically.

1.2. There is a entire curve $h:\C\to X$ with Zariski-dense image.

1.3. There is a entire curve $h:\C\to X$ with metrically dense image.

1.4. Any two generic points of $X$ are joined by some entire curve.

1.5. Any countable subset of $X$ is contained in the image of some entire curve.

\item
An arithmetic analogue is the following:

$X$ is potentially dense (i.e.: $X(k)$ is Zariski-dense for $k$ a sufficiently large number field $k$ of definition of $X)$.

\item
  Moreover we conjecture that
  $\pi_1(X)$ is almost Abelian if $X$ is special (\cite{Ca04}, Conjecture 7.1).
  \end{enumerate}

\end{Conjecture}

$\bullet$ The property 1.4 means that special manifolds are expected to be analogues of rationally connected manifolds, replacing algebraic maps from $\C\to X$ by their transcendental version: entire curves $h:\C \to X$, with $h$ holomorphic, non-constant.

$\bullet$ In general, we have the obvious implications:

1.5 $\Longrightarrow$ 1.3 $\Longrightarrow$ 1.1, 1.3 $\Longrightarrow$ 1.2, and
1.5. $\Longrightarrow$ 1.4 $\Longrightarrow$ 1.1.

$\bullet$ If $X$ is unirational, all the properties above ($1.1-1.5, 2$ and $3$)
are satisfied.

$\bullet$ 
If there exists a surjective holomorphic map
$p:\C^n\to X$, then all these properties 1.1-1.5 are satisfied.

$\bullet$ Properties 1.1 and 1.4
are obvious for rationally connected
manifolds.

The main purpose of this article is to show that every rationally connected
manifold admits a dense entire curve.
In fact, we prove  1.2, 1.3 and 1.5 for rationally connected manifolds.
(Theorem~\ref{rc-main}).
Thus all the
properties
$1.k$ ($1\le k\le 5$) hold for every rationally connected manifold, i.e.,
we confirm part $1.5$ of the conjecture for the simplest class
of special manifolds,
namely rationally connected manifolds.

$\bullet$ It is proved in \cite{BL} that a projective surface is 'special' in the above sense if and only if it is $\C^2$-dominable (With the possible exception of non-elliptic and non-Kummer $K3$ surfaces.). This might however be a low-dimensional phenomenon, and it is not expected to remain true in dimension $3$.

\begin{remark} Vojta
\cite{Vj}
  has introduced a dictionary between entire curves and infinite sequences of rational points. The conjectures above imply that if $X$ is defined over a number field, the set of $k'$-rational points is infinite for some finite extension $k'/k$ if and only if $X$ contains an entire curve, Vojta's analogy thus becomes an equivalence between entire curves and arithmetic geometry. The Nevanlinna version of the Hilbert property discussed in \S9 is another illustration of this statement.
\end{remark}

\section{Entire Curves in Rationally Connected Manifolds}

We recall here classical results used in the next subsections.

\subsection{Comb smoothing}

  We recall the technique of {\em ``comb smoothing''},
  introduced in \cite{KMM}, see also \cite{D}.

 In the sequel, we always see $\P_1=\C\cup \{ \infty\}$, ie, with a distinguished point at infinity.
 
 Let $C$ be the reducible curve given by glueing two rational curves together
 transversally in one point.

 Let $V$ be the variety obtained from blowing up the point $(\infty,0)$
 in $\P_1\times\C$. Using $[y_0:y_1],t$ as coordinates
 on $\P_1\times\C$, we
 may embedd $C$ into $V$ as the total transform of $t=0$.
 With $L$ denoting an affine line, the function $t$ realizes $V$ as a
 flat projective $L$-scheme,
 with a special fibre isomorphic to $C$ and all
 other fibers isomorphic to smooth rational curves.
 By abuse of notation, we denote this special fiber again as $C$.

 In explicit coordinates
 \begin{align*}
   V &=\{([y_0:y_1];t;[s_0:s_1])\in\P_1\times\C\times\P_1
   :s_1y_0=s_0y_1t\}
   \\
   &\mapsto t\in L
 \end{align*}
 For a given morphism $f:C\to X$ with values in a projective manifold $X$,
 the morphism $f$ can be extended to a holomorphic map $F$ defined on an open
 neighborhood of $C$ in $V$ if the following condition $(*)$ is fulfilled:

 \medskip
 
\noindent\fbox{
 $f^*TX$ is spanned by global sections on $C$ and $H^1(C,f^*TX)=\{0\}$
}

\medskip

Let $r\geq 0$ be an integer, recall (\cite{D}, Definition 4.5) that a
(pa\-ra\-me\-trised) rational curve $f:\P_1\to X$ is said to be $r$-free if $f^*(TX)\otimes \mathcal O_{\P_1}(-(r-1))$ is ample on $\P_1$, or equivalently, if $H^1(\P_1,f^*(TX)\otimes \mathcal O_{\P_1}(-(r+1)))=0$. This is an open condition on the space of such maps.

The condition ``$1$-free'' is frequently called ``very free''.
 
   If $f:C\to X$ is a morphism whose restriction to the two irreducible
   components of $C$ are $r_1$-free resp.~$r_2$-free,
   then the condition $(*)$ is fulfilled.
   Hence $f$ extends to a holomorphic map $F$ defined on some neighborhood
   of $C$ in $V$. For each fiber of $V\to L$ contained in this neighborhood
   we obtain a rational curve $f_t$.
   Moreover in this case the curves $f_t$ are $(r_1+r_2)$-free rational
   curves.

   Summarizing, we have the following:
   \begin{proposition}\label{comb-smoothing}
     Let $X$ be a projective manifold, equipped with a distance function $d$
     induced by a hermitian metric on $X$, and let  $g_i:\P_1=\C\cup\{\infty\}\to X$ be
     $r_i$-free rational curves ($i=1,2$, $r_i\ge 0$) with $g_1(\infty)=g_2(0)$.

     Then for every $R,R',\epsilon>0$ there is a $(r_1+r_2)$-free rational curve
     $c:\P_1=\C\cup\{\infty\}\to X$ and a parameter $\lambda>0$ such that
     $d(c(z),g_1(z))<\epsilon$ for all $z$ with $|z|<R$ and
     $d(c(\lambda z),g_2(z))<\epsilon$ for all $z$ with $|z|>R'$.

     Moreover given $p_0,p_1,\ldots p_{r_1}\in \C$ and $q_1,\ldots,q_{r_2}\in\C^*$
     the curve
     $c$ may be chosen in such a way that there are points
     $\tilde p_i,\tilde q_i$ on $\C$ with
     $c(\tilde p_i)=g_1(p_i)$, $|\tilde p_i-p_i|<\epsilon$
     and $c(\tilde q_i)=g_2(q_i)$.
   \end{proposition}

   \begin{proof}
     Let $V$ be as above. We embed the first rational curve into $V$
     as
     \[
     \C\cup\{\infty\}\ni x \mapsto ([1:x],0,[1:0])
     \]
     and the second one as
     \[
     \C\cup\{\infty\}\ni x \mapsto ([0:1],0,[1:x]).
     \]
     Now the pair $(g_1,g_2)$ defines a morphism from $C$ to $X$ which
     extends to holomorphic map $F$ defined on  an open
     neighborhood. Thus for small enough $t$
     we obtain a rational curve $f_t:\C\cup\{\infty\}\to X$
     as
     \[
     \C\cup\{\infty\}\ni x \mapsto
     F\left( [1:x],t,[1:xt]\right).
     \]
     We define $\lambda=\frac 1t$ and obtain
     \[
     \C\cup\{\infty\}\ni \lambda x \mapsto
     F\left( [1:\lambda x],t,[1:\lambda xt]\right)=
     F\left( [t: x],t,[1:x]\right).
     \]
     It is now easy to check explicitly that for sufficiently small $t$ the
     rational curve $f_t$ has the desired properties.
   \end{proof}

\subsection{Dense entire curves in Rationally connected manifolds}

We present our
main result on the existence of dense entire curves
in rationally connected manifolds. Later we deduce
from it stronger forms by applying it to projective jet bundles of $X$ (see corollary \ref{rc-corollary1} below).

\begin{theorem}\label{rc-main}
  Let $X$ be a rationally connected complex projective manifold and
  let $A$ be a closed analytic subset of codimension at least two,
  and let $M$ be a countable subset of $X\setminus A$.

  Then there exists an entire curve $h:\C\to X$ such that
  \begin{enumerate}
 
    \item
      $M\subset h(\C)$,
    \item
    $h(\C)\subset X\setminus A$,
  \end{enumerate}
        In particular, if $M$ is dense in $X$, so is $h(\C)$.
\end{theorem}

\begin{corollary}
  Rationally connected manifolds satisfy all the properties of part $(i)$
  of Conjecture~\ref{conjspec}.
\end{corollary}

\begin{assumptions}\label{ass-1}
In the sequel, $X$ always denotes a rationally connected complex projective
manifold, $A$ a closed analytic subset of $X$ of codimension at least $2$,
and $d(\ ,\ )$ denotes a distance function on $X$ induced by some smooth hermitian
metric on $X$.
\end{assumptions}
 
\begin{lemma}\label{rc-lemma} Let $X, A, d$ be as defined above
  and $r\in\N$.
  Let  $f:\P_1\to X$
  be a $(r-1)$-free rational curve such that $f(\P_1)\cap A=\emptyset$.
  Let $a_i\in\C, i=1,\dots,r$ be given, let $R>0,\varepsilon>0$ be given, with $R\geq \vert a_i\vert+\varepsilon, \forall i$. Let $q\in X\setminus A$.
 
 There then exists a rational curve $h:\P_1\to X$ and $a'_i\in \C, i=1,\dots, r+1$ such that:
 
1. $h$ is $r$-free.

2. $h(\P_1)\cap A=\emptyset$.

3. $d(h(z),f(z))\leq \varepsilon$ if $\vert z\vert \leq R$.

4. $h(a'_i)=f(a_i), i=1,\dots, r$.

5. $\vert a'_i-a_i\vert \leq \varepsilon,i=1,\dots, r$.

6. $h(a'_{r+1})=q$.

 \end{lemma}

\begin{proof} We first prove the existence of a $1$-free curve $g:\P_1\to X$ avoiding $A$ and containing $q$ and $f(\infty)$. Indeed, from \cite{D}, Corollary 4.28, we get the existence of an $1$-free curve $g:\P_1\to X$ such that $g(0)=f(\infty)$ and $g(1)=q$.
  Since the deformations of $g$ are unobstructed, there is a holomorphic map $G:T\times \P_1\to X$, where $T$ is an open neighborhood of $0$ in $H^0(\P_1,f^*(TX)\otimes \mathcal O_{\P_1}(-\{0,1\}))$ inducing deformations $G(t,.)=g_t$ of $g$ fixing $g(0)$ and $g(1)$, with $dG(0, z)(t)=t(z)$, for $z\in \P_1$, and $t\in H^0(\P_1,f^*(TX))$.
  Now $g^*(TX)$ is generated outside of $0$ and $1$ by its sections, and so $G$ has maximal rank near $\{0\}\times (\P_1\setminus\{0,1\})$, so that $G^{-1}(A)$ has codimension at least $2$ in $T\times \P_1$.
  Its projection in $T$ has thus codimension at least one, and hence the generic $g_t$ avoids $A$.

   We now apply the``comb smoothing'' technique
   (see Proposition~\ref{comb-smoothing}.)  to the comb
   defined by the rational curves $f$ and $g$. We obtain
   a $r$-free rational curve $h$ fulfilling the properties
   1,3,4,5,6.
   Since $f(\P_1)\cup g(\P_1)$ do not intersect $A$, we see that
   $C\cap F^{-1}(A)$ (in the notation of prop.~\ref{comb-smoothing})
   is empty.
   It follows that $F^{-1}(A)$ does not intersect an open neighborhood
   of $C$ in $V$. Therefore we may choose $h$ such that
   $h(\P_1)\subset X\setminus A$
   (i.e., such that condition 2 is fulfilled).
  \end{proof}
  
\begin{proof} We now prove Theorem \ref{rc-main}. Let $M:=\{x_1,x_2,\dots,x_n,\dots\}$ be the given sequence in $X\setminus A$. We  construct inductively a sequence of rational curves
  $h_n:\P_1=\C\cup\{\infty\}\to X$, $R_n>0$, $\varepsilon_n>0$, $a_n^k\in \C,\forall k\geq n\geq 1$ such that:

  0. $\varepsilon_n\to 0$ when $n\to +\infty$, the sequence $\varepsilon_n$ being decreasing.
  
  1. $h_n$ is $n$-free.
  
  2. $d(h_n(\P_1),A)\geq \varepsilon_n>0$, in particular: $h_n(\P_1)\cap A=\emptyset$.
  
  3. $d(h_{n+1}(z),h_n(z))< \frac{\varepsilon _n}{2^{n+1}}$ for $\vert z\vert \leq R_n$.
  
  4. $R_{n+1}\geq R_n+1$.
  
   5. $h_k(a_n^k)=x_n,\forall k\geq n\geq 1.$
   
    6. $\vert a_n^{k+1}-a_n^k\vert\leq \frac{\varepsilon _n}{2^{k+1}},\forall k\geq n\geq 1$.

  7. $\vert a_n^n\vert \leq R_n$

  Assume these data are given. From conditions 0, 3 and 4, we see that the maps $h_n$ converge uniformly on compact subsets of $\C$ to a holomorphic map $h:\C\to X$.
  For a fixed number $n$, the sequence $a_n^k$ is a  Cauchy sequence
  by condition 6, hence convergent.
  We define $a_n:=\lim_{k\to +\infty}a_n^k$.
  Condition 5 implies that $h(a_n)=x_n,\forall n$.
  Thus $M\subset h(\C)$.

   From condition 3 we obtain that $d(h(z),h_n(z))\le \frac{\epsilon_n}{2^n}$
   for all $z\in\C$ with $|z|\le R_n$.
   In combination with condition 2 this yields
   $d(h(z),A)\ge \epsilon_n\left( 1 - 2^{-n}\right)>0$
   for all $z\in\C$ with $|z|\le R_n$.
   Since this holds for all $n$, we obtain $d(h(z),A)>0$ for all
   $z\in\C$, i.e.,
   $h(\C)\cap A=\emptyset$.

  We now construct the sequences $h_n,R_n,\varepsilon _n, a_n^k$ enjoying the above properties 1-7 by induction on $n$.
  This will prove the Theorem \ref{rc-main}.
  
  For $n=1$, just take for $h_1$ any $1$-free rational curve avoiding $A$, and with $h_1(1)=x_1$, and choose $a^1_1=R_1=1$. Choose $\varepsilon_1=\frac{d(h_1(\P_1),A)}{2}$.

  Fix $n$.
  Assume $h_m,\varepsilon_m,R_m,a_m^k$ already constructed for all $k,m$ with
  $m\leq k\leq n$, and satisfying the properties 1-7 above.
  We apply the Lemma \ref{rc-lemma} with $f=h_n$, $a_i=a_i^n$, 
  $R=R_n$ and $\epsilon=\frac{\epsilon_n}{2^{n+1}}$.
  This yields an $(n+1)$-free rational curve $h_{n+1}=h$.
  We set $a_k^{n+1}=a_k'$ (in the notation of lemma~\ref{rc-lemma}).
  We finally choose $R_{n+1}>\max\{R_n+1,|a_{n+1}^{n+1}|\}$ and
  $\varepsilon_{n+1}<\min\{\frac{\varepsilon_n}{2^{n+1}},d(h_{n+1}(\P_1),A)\}$.
\end{proof}

  A slight modification of the proof gives:

\begin{corollary}\label{cor-dense-ini}
  Let $X$ be a rationally connected projective manifold with a free rational
  curve $c:\P_1\to X$ and a distance function $d(\ ,\ )$ induced by some
  hermitian metric.

  Then for every $R>0,\epsilon>0$ there exists an entire curve $h:\C\to X$
  with dense image and such that $d(h(z),c(z))<\epsilon$ for all $z\in\C$
  with $|z|<R$.
\end{corollary}

\begin{corollary}\label{rc-corollary1} Let $X,A,M$ be as above
  in \ref{ass-1}, let $D\subset X$ be a reduced divisor (not necessarily irreducible nor connected), and let $m>0$ be an integer.
  For each $x_n\in M$, fix an (unparametrized) $m$-jet $j_n$ of map from the unit disc in $\C$ to $X$ at $x$. The map $h$ of Theorem \ref{rc-main} which sends $a_n\in \C$ to $x_n$ can be choosen to have its $m$-jet at $a_n$ coincide with $j_n$, this for each $n$, and moreover: to avoid $D^{sing}$, and to be transversal to $D^{reg}$ at each point where $h(\C)$ meets $D^{reg}$. (Here $D^{sing}$ is the singular locus of $D$, and $D^{reg}$ is its regular part).
  \end{corollary}
  
  \begin{proof} We prove this first for $m=1$: let $\pi:X_1:=\P(TX)$ be the natural projection. Since $X$ is rationally connected, so is $X_1$.
  
  Lift arbitrarily the sequence $M$ to the given sequence $M_1$ in $X_1$, and put $A_1:=[\pi^{-1}(A\cup D^{reg})]\cup \P(TD^{reg})\subset X_1$. Apply Theorem \ref{rc-main} to $X_1,A_1,M_1$ to get the result, since $A_1$ is algebraic, closed, and of codimension at least $2$ in $X_1$. Then proceed inductively on $m$, by defining $X_{m}:=(X_{m-1})_1$ to get the general case if $m>1$.
  \end{proof}

\begin{remark} By composing $h$ with a suitable Weierstrass product $w:\C\to \C$, one can even realise the jets $j_m$ themselves, and not only at the level of their projectivisation. This gives an entire curve analogue of the Weak Approximation Property in arithmetic geometry on rationally connected manifolds. This connection with the Nevanlinna Hilbert Property treated below was pointed to us by P. Corvaja.
\end{remark}

  One can refine Corollary~\ref{rc-corollary1}, by imposing to $h$ to meet $D$ transversally at every of its regular intersection point:

  \begin{theorem}\label{rc-main-2}
    In the set-up of theorem~\ref{rc-main}, assume that in addition a
    (reduced)
    hypersurface $D$ in $X$ is given. Then we may chose the entire
    curve $h:\C\to X$ provided by theorem~\ref{rc-main} in such a
    way that {\em every} intersection point of $h(\C)$ and $D$
    is transversal.
  \end{theorem}

Before proving the theorem, we need some auxiliary lemmata.
  
  \begin{lemma}\label{5.8}
    Let $D$ be a reduced hypersurface in a projective manifold $X$ and
    let $u:\P_1\to X$ be a $2$-free rational curve.

    Then there exist arbitrarily small deformations $\tilde u$ of $u$
    such that
    every intersection point of $\tilde u(\P_1)\cap D$ is transversal.
  \end{lemma}

  \begin{proof}
    A  generic deformation of $u$ will avoid
    any given analytic subset of codimension at least two.
    Therefore a generic deformation of $u$ avoids $\Sing(D)$.
    Now let $p\in\P_1$ with $u(p)\in D\setminus \Sing(D)$.
    By assumption $u^*TX\tensor I(2\{p\})$ is an ample vector bundle
    on $\P_1$. By the theorem of Grothendieck $u^*TX$ is a direct sum
    of line bundles. We choose a direct summand $L$ of $u^*TX$
    which is complimentary to $u_p^*TD$ at $p$.
    Now we take a section of $u^*TX$ which vanishes at $p$, but has
    non-zero derivative in the direction of $L$. Then the corresponding
    deformations yield rational curves which intersect $D$ with multiplicity
    $1$ in $p$. Since $D\cap u(\P_1)$ is finite
    (in fact its cardinality is bounded by $\deg(u^*{\mathcal L}(D))$),
    it follows that a generic deformation will have transversality at every
    point of intersection.
  \end{proof}

  \begin{lemma}
    Let $X$ be a complex manifold, $D$ a (reduced) hypersurface on $X$,
    and $h:\C\to X$ an entire curve.
    Let $G$ be a relatively compact open subset of $\C$ with $h(\partial G)
    \cap D=\emptyset$. Let $d$ be a distance function on $X$ induced by
    some hermitian metric.

    Then there exists a number $\delta>0$ such that the following
    two properties hold for every entire curve $\tilde h:\C\to X$
    with $d(\tilde h(z),h(z))<\delta$, $\forall z\in\bar G$:
    \begin{enumerate}
    \item
      $\tilde h(\partial G)\cap D=\emptyset$.
    \item
      $\deg_G(h^*D)=\deg_G(\tilde h^*D)$.
    \end{enumerate}
  \end{lemma}

  \begin{proof}
    Follows easily using the theorem of Rouch\'e.
  \end{proof}

  \begin{lemma}\label{5.10}
    Let $X$ be a complex manifold, $D$ a (reduced) hypersurface on $X$,
    and $h:\C\to X$ an entire curve.
    Let $G$ be a relatively compact open subset of $\C$ with $h(\partial G)
    \cap D=\emptyset$. Let $d$ be a distance function on $X$ induced by
    some hermitian metric. Assume that $h^*D$ is reduced (i.e.~all the
    multiplicities are one).

    Then there exists a number $\delta>0$ such that the following
    two properties hold for every entire curve $\tilde h:\C\to X$
    with $d(\tilde h(z),h(z))<\delta$, $\forall z\in\bar G$:
    \begin{enumerate}
    \item
      $\tilde h(\partial G)\cap D=\emptyset$.
    \item
      $\tilde h^*D$ is reduced on $G$.
    \end{enumerate}
  \end{lemma}
  \begin{proof}
    Let $p_1,\ldots,p_s$ be the points of $h^*D$.
    We choose small balls $B_i$ such that $p_i\in B_i\rc G$
    and such that $p_j\not\in\overline{B_i}$ for $i\ne j$.
    Then we choose $\epsilon_i>1$ such that $\deg_{B_i}(\tilde h^*D)=1$
    for every entire curve $\tilde h:\C\to X$ with
    $d(\tilde h(z),h(z))<\epsilon_i$, $\forall z\in B_i$.
    (This we can do due to the preceding lemma).
    Finally we choose $\epsilon=\min_i\epsilon_i$.
  \end{proof}

  \begin{proof}[Proof of theorem~\ref{rc-main-2}]
    We proceed as in the proof of theorem~\ref{rc-main} with the
    following modifications. We choose the rational curve $h_n$ such that
    $h_n^*D$ is reduced on $\{z\in\C:|z|\le R_n\}$
    (This is possible due to Lemma~\ref{5.8}). Thanks to
    Lemma~\ref{5.10} there exists a
    number $\delta_n>0$ such that any entire curve $\tilde h:\C\to X$
    with $d(h_n(z),\tilde h(z))<\delta_n$, $\forall |z|\le R_n$ preserve this
    property (i.e., $\tilde h^*D$ is reduced on the disc with radius
    $R_n$). We choose $\epsilon_n$ such that $\epsilon_n<\delta_n$.
    In this way we finally arrive at an entire curve $h$ with
    $h^*D$ being reduced everywhere.
  \end{proof}

 \begin{remark} The arithmetic version (i.e.: the potential density of Rationally connected manifolds $X$ defined over a number field) of the Theorem
   \ref{rc-main} is open when $X$ is not known to be unirational, even for conic bundles over $\P_n,n\geq 2$. There are (at least) two cases known: smooth quartics in $\P_4$ (\cite{HT}), and some conic bundles over $\P_2$ (\cite{BMS}, where $X(\Q)$ is shown to be already Zariski-dense).
\end{remark}

  \subsection{An orbifold situation.}

In some situations, we can impose tangency conditions as well. These are simple instances of an extension of the preceding results to the `orbifold pair' situation, for which the appropriate techniques have not been presently developped.

 \begin{corollary}\label{rc-corollary 3} Let $S\subset \P_3$ be a smooth sextic surface. There exists entire curves $g:\C\to \P_3$ with metrically dense image $C$, and such that $C$ is tangent to $S$ at each point where it meets $S$. 
 \end{corollary}
 
 \begin{proof} Let $\pi:X\to \P_3$ be the double cover
   ramified along $S$:
   $X$ is rationally connected since it is Fano. Let $h:\C\to X$ be a dense entire curve. Then $g:=\pi\circ h:\C\to \P_3$ is a dense entire curve with the claimed property.\end{proof}
 
 We don't know how to prove Corollary \ref{rc-corollary 3} directly on $(\P_3,S)$, without going to the double cover.

 \begin{remark} Observe that
   in the situation of the above corollary
   $K_{\P_3}+(1-\frac{1}{2}).S$ has negative degree (i.e.: the smooth `orbifold pair' $(\P_3,(1-\frac{1}{2}).S)$ is Fano). The statement of the preceding corollary thus means that this orbifold pair admits metrically dense `orbifold entire curves' (in the `divisible' sense). This is conjectured in \cite{Ca} to be the case for arbitrary smooth orbifolds which are Fano (or, more generally, Rationally connected in the orbifold sense).
 \end{remark}

 \subsection{Brody curves}

 An entire curve $f:\C\to X$ is a ``Brody curve'' if its derivative
 is uniformly bounded with respect to the euclidean metric on $\C$ and an
 arbitrary hermitian metric on $X$. The theorem of Brody states that
 on a compact complex manifold $X$ there is a
 non-constant Brody curve if and only if there is a non-constant entire
 curve.

 However, there is no such statement concerning {\em dense} curves.
 Indeed, an Abelian threefold blown-up along a curve is constructed in \cite{W}, 
 which admits (lots of) dense entire curves, although every Brody curve
 is contained in the exceptional divisor. The existence of dense Brody curves is thus not a bimeromorphic property.
  
 Our results give no information about the existence of (Zariski-) dense Brody
 curves on arbitrary rationally connected manifolds (although they obviously exist on some of them).

 \section{Normal rational surfaces.}

 If $X$ is a singular complex algebraic variety with desingularization
 $\pi:\tilde{X} \to X$, then there exists a (Zariski-)dense entire
 curve $h:\C\to X$ if and only if there exists a (Zariski-)dense
 entire curve $\tilde h:\C\to \tilde X$. (This is clear, because we can
 lift every
 holomorphic map $h:\C\to X$ whose image
 $h(\C)$ is not contained in the singular locus $X^{sing}$.)

 However, it is a much more delicate question whether there exists
 a (Zariski-)dense entire curve $h:\C\to X$ which avoids the
 singular locus $X^{sing}$ of $X$.

 In this section we prove the existence of such curves for certain
 normal rational surfaces (Theorem~\ref{rsurf}). Note that most normal rational surfaces have no Zariski dense
 entire curves in their smooth locus. See examples \ref{ueno} and ~\ref{rational-ex} at the end of this
 section. 
 
 We consider the following situation: $S$ is a normal projective surface with only quotient singularities, such that $-K_S=\Delta$, an effective nonzero $\Q$-divisor on $S$ with $\coeff(\Delta)\leq 1$, where $\coeff(\Delta)$ is the largest of the coefficients of the irreducible components of the support of $\Delta$. 
    
  This permits to apply to $S$ the MMP, which preserves the condition on $\coeff(\Delta)$, and abuts, after finitely many $K$-negative contractions, to a pair $(S',\Delta')$ such that either $K_{S'}$ is nef, or anti-ample with Picard number $1$, or  admitting a Fano-contraction $f:S'\to B$ on a smooth projective curve $B$ with relative Picard number $1$.
 This allows us to characterise the existence of dense entire curves in the
  smooth locus of $S$.
  
  \begin{theorem}\label{rsurf}
    Let $S$ be a normal rational surface with only quotient singularities, and $F\subset S$ be a finite subset containing the singular points of $S$. Let
    $\Delta$ be an effective nonzero
    $\Q$-divisor such that  $(K_S+\Delta)$ is $\Q$-trivial.
    \begin{enumerate}
      \item
        If $\coeff(\Delta)<1$,
        then
        for any two generic points of $S\setminus F$,
        there is a $1$-free rational curve through these two points
        avoiding $F$.

      \item
        If $\coeff(\Delta)\leq 1$,
      then there is an entire curve $h:\C\to S\setminus F$ with dense image.
    \end{enumerate}
  \end{theorem}

  \begin{proof}
   Applying the Minimal Model Program (MMP) for surfaces
    (see \cite{KM'}, Theorem 3.47, 4.11 and 4.18), we obtain a birational map $\mu:S\to S'$ to a normal surface
    $S'$ with only quotient singularities
    and a $\Q$-divisor $\Delta':=\mu_*(\Delta)$ on $S'$ with $\coeff(\Delta')\leq 1$,
  and $K_{S'}+\Delta'$ $\Q$-trivial such that one of
    the following possibilities hold:

    \begin{enumerate}
  \item
  The canonical divisor $K_{S'}$ is nef.
  \item
    The anticanonical divisor $-K_{S'}$ is ample.
    
      \item
    There is a Fano fibration $f:S'\to B=\P_1$ of relative Picard number $1$.
    \end{enumerate}
    In our situation, the first case does not occur.
    Indeed, in this case we had: $0\leq \kappa(S')=
    \kappa(S)=\kappa(S,-\Delta)=-\infty$,
    since: $\Delta$ is $\Q$-effective and nonzero, the MMP preserves the Kodaira dimension, and $\kappa(S')\geq 0$ if $K_{S'}$ is nef.

     In the second case, our first claim has been proved in \cite{Xu}.
     Claim 2 now follows from claim 1 just as in the smooth case,
     proved in Theorem \ref{rc-main}. 
     
     In the third case the
     existence of a dense entire curve and of rational curves with the said properties follow from
     Proposition~\ref{fano-rational}.
 \end{proof}

  \begin{remark} When $\Delta=0$ above, one also expects the existence of dense entire curves, but this is an open question. See  Example \ref{ueno} for an example in which this conclusion does not follow from the arguments given here, but from an explicit description. 
  
  It were interesting also to get a criterion to recognise from the initial data $(S,\Delta)$ itself, whether one abuts, for some suitable MMP, to the $\Bbb Q$-Fano case with $\rho=1$, or to the fibered case for any MMP-sequence of contractions. \end{remark}

  We need some preparatory lemmas.
  
  \smallskip
  
    {\bf Assumptions 1.}
    
    \smallskip

 $S$ denotes a projective normal surface with only quotient
   singularities, and there exists a surjective Fano-fibration $f:S\to B$
    of relative Picard number $1$. 
This implies that
   all fibers $S_b=f^{-1}\{b\}$ are irreducible.
   The $\Q$-divisor $\Delta_f:=\sum_b(1-\frac{1}{m_b})\{b\}$
     (with $m_b$ being the multiplicity of the fiber $f^{-1}(b)$
     for each $b\in B$)
     denotes the `orbifold base' of $f$.

   \begin{lemma}\label{cw}
     Let $S$ be a normal surface with a holomorphic map $f:S\to B$
     with irreducible fibers. If there exists an entire curve $h:\C\to S^{reg}$ such that $f\circ h:\C\to B$ is nonconstant, then $deg(K_B+\Delta_f)\leq 0$.
    Thus: either $B$ is elliptic and $\Delta_f=0$, or $(B,\Delta_f)=(\P_1,\Delta_f)$ is in the short classical `platonic' list of orbifolds on $\P_1$ with
     (semi-)negative canonical bundle.
   \end{lemma}

   \begin{proof}
     This follows from \cite{CW}, since $f\circ h$ is then an
     orbifold morphism to $(B,\Delta_f)$.
      \end{proof}

        We shall next use the following properties:
   
   1. If $m_b=1$, $S_b\cong \P_1$ is smooth (\cite{KM}, 11.5.1), and $f$ is thus a $\P_1$-bundle over $B$ if and only if $\Delta_f=0$.
   
   2. Assume\footnote{This is true unless $B=\P_1$, and the support of $\Delta_f$ consists of one or two points with different multiplicities.
     However, in these two cases
     there is still  such a dominant finite algebraic map
     $\beta:\C\to B$ which suffices for our purposes.} there exists a finite ramified cover $\beta:B'\to B$ which ramifies at order exactly $m_b$ over $b$, this for each $b\in B$. Let $\sigma:S'\to S$ and $f':S'\to B'$ be deduced by base-changing $f$ by $\beta$, and normalising the base-change $S\times_BB'$. Then $f'$ has generically reduced fibres, $S'$ has still quotient singularities, with $K_{S'}+\Delta'$ $\Q$-trivial  (\cite{KM'}, Proposition 5.20 and Proposition 4.18), and $\coeff(\Delta')\leq 1$. We can still apply to $S'$ the MMP and abut now to a fibration $f":S"\to B'$ which has only reduced fibres, and is thus a $\P_1$-bundle.
   
   \begin{lemma}\label{base change} Let $f:S\to B$ be as in the assumptions above. Let $f':S'\to B'$ be deduced from $f$ by a finite base-change $B'\to B$. Assume that all fibres of $f':S'\to B'$ have a reduced component. Let $F\subset S^{reg}$ be a finite set. Then $S^{reg}\setminus F$ contains a dense entire curve (resp. a $1$-free rational curve through any two of its generic points) if $B'$ is elliptic (resp. rational).
   \end{lemma}
   
\begin{proof} It is sufficient to show this for $(S')^{reg}\setminus F'$, where $F'$ contains the inverse image in $S'$ of $S^{sing}$ and of $F$. But $S'$ dominates a $\P_1$-bundle over $B'$, after the preceding observations, which implies the claims.\end{proof}

The following corollary then follows from the existence of a suitable `orbifold-\'etale' base change $B'\to B$, with $B'$ either elliptic or rational, when $(B,\Delta_f)$ is in this list of `platonic' orbifolds.

\begin{corollary}\label{surf'} Let $f:S\to B$ be a Fano fibration of relative Picard rank one, $S$ a normal projective surface with only quotient singularities. The following properties are then equivalent:

  1. $\deg(K_B+\Delta_f)<0$ (resp. $\deg(K_B+\Delta_f)\leq 0)$, where $(B,\Delta_f)$ is the orbifold
         basis of $f$.

          2. There exists a rational curve (resp. an entire curve) on $S^{reg}$ which is not $f$-vertical.\footnote{i.e.: not contained in some fibre of $f$.}        
         
          3. There exists a $1$-free rational curve through any
         two generic points of $S$ (resp. a dense entire curve)
         which is contained in $S^{reg}$.
\end{corollary}

\medskip

$\bullet$ We shall now relate $K_S$ to $K_B+\Delta_f$ by means of a `canonical bundle' formula for $K_S$.

     \medskip

     Let $r:\tilde{S}\to S$ be a minimal resolution, and $
     t:\tilde{S}\to S_0$ be a relative minimal model of $\tilde{S}$ over $B$, so that $f_0:S_0\to B$ is a $\P_1$-bundle, with $f_0\circ t=f\circ r$.
     Then $S_0=\Bbb F_m$, the Hirzebruch surface of index $m$
     for some $m\in\N_0$.
     In this case, we denote by $G$ (resp. $G')$ the section of $f_0$ with $G.G=-m$ (resp. a section of $f_0$ with $G.G'=0$, or equivalently, disjoint from $G)$.
     We have: $K_{S_0}=-(G+G')+(f_0)^*(K_B)$. We shall then denote with $C$ (resp. $C')$ the strict transforms in $S$ of $G,G'$.

     \begin{lemma}\label{canonical bundle formula}
       Under the above assumptions,
       we have: 
       
       $K_S\sim f^*(K_B+\Delta_f)-(C+C')$.
       Moreover, $C$ and $C'$ are still disjoint.
   \end{lemma}
   
   \begin{proof} Each singular fibre of the fibration $f\circ r:\tilde{S}\to B$ has two reduced components, each of which is met once by either $\tilde{G}$ or $\tilde{G'}$, the strict transforms of $G,G'$ in $\tilde{S}$. The contraction $r:\tilde{S}\to S$ then sends each of these components to $2$ distinct singular points of the corresponding fibre of $f$, by \cite{KM}, 11.5. 5, which implies that $C,C'$ are disjoint, too. The first claim then follows from the fact that $F_b:=(S_b)_{red}$ is irreducible for each $b\in B$, so that $K_S=a.F_b-C-C'$ locally near $F_b$, with $a=(m_b-1)$ by adjunction, since $F_b$ has multiplicity $m_b$ in $f$. \end{proof}

   \begin{lemma}\label{intersections}
     Under the above assumptions $C.C'=0$ and  $C'\sim C+k\Phi$
     for some $k \geq 0$ (possibly exchanging $C,C')$.
     Moreover, for each $f$-horizontal irreducible curve $H\subset S$ different from $C,C'$, we have: $H\sim d(C'+\ell \Phi)$ with $d:=H.\Phi>0$, $\ell=\frac{H.C}{d}\geq 0$.
   \end{lemma}

   \begin{proof}
     $S$ is $\Q$-factorial (because it has only quotient singularities)
     and its Picard number is $2$, since $f:S\to B$
     has relative Picard number one.
     Hence every divisor is numerically equivalent to a linear
     combination of $C'$ and $\Phi$. The assertions of
     the lemma are now obvious.
   \end{proof}

   {\bf Assumptions 2.}
   
   Let $f:S\to B$ be a Fano fibration with $S$ normal and having only quotient singularities. Assume now that $K_S+\Delta$ is $\Q$-trivial, where $\Delta$ is a nonzero effective $\Q$-divisor on $S$. We write $\Delta=\Delta^h+\Delta^v$, where $\Delta^h$ (resp. $\Delta^v)$ is the $f$-horizontal (resp. $f$-vertical) part of $\Delta$.
   
   We write $\coeff^h(\Delta)$ for the largest coefficient in $\Delta$ of the irreducible components of $\Delta^h$. 
   
   \begin{lemma}\label{K_S+D} Under the above assumptions,
     $\deg(K_B+\Delta_f)\le 0$.
     Moreover, $\deg(K_B+\Delta_f)<0$ if $coeff^h(\Delta)<1$.
   \end{lemma}
   
   \begin{proof}
    Since: $\Delta=\Delta^h+\Delta^v \sim K_S\sim (C+C')-f^*(K_B+\Delta_f),$
     $\sim$ denoting linear equivalence, we have: 
     $-f^*(K_B+\Delta_f) \sim \Delta^v + \Delta^h-(C+C')$,
     and it is sufficient to show that:
      $\Delta^h-(C+C')$ is $\Q$-effective.
     We may write: $\Delta^h=aC+a'C'+R$, where $R=\sum_ja_j\frac{C_j}{d_j}$
     with $C_j$ being irreducible multisections of $f$ distinct from $C,C'$,
     each of degree $d_j>0$ over $B$, and $a_j>0$ are rational numbers.
     
     We have, for each $j$: $C_j=d_j(C'+k_j\Phi)$ for some $k_j\in\Q$.
     Thus $R=\sum_j a_j(C'+k_j\Phi)$. Calculating intersection numbers with
     $\Phi$ gives:
     $2=(C+C').\Phi=\Delta^h.\Phi=a+a'+\sum_j a_j$.
     Therefore: 
     
     $ \Delta^h-(C+C')\sim aC +a'C'+(2-a-a')\sum_j(C'+k_j\Phi)-(C+C')$
     
     $=(1-a)(C'-C)+(2-a-a')\left(\sum_j k_j\right)\Phi.$

     Recall that  $k_j\ge 0$.
     Furthermore $(C'-C)$ is effective (lemma~\ref{intersections}).
     Hence the above equation implies that
     $(\Delta^h-C-C')$ is $\Q$-effective, and so is $(\Delta-C-C')$
     which implies: $\deg(\Delta_f+K_B)\le 0$.

     If $\coeff(\Delta)<1$, we have: $a,a'< 1$.
     Consequently $(1-a)>0$ which implies that $(\Delta^h-C-C')$
     is non-zero, and so: $\deg(\Delta_f+K_B)<0$.
     \end{proof}

   We have the following criterion for the existence of
   $1$-free rational curves (resp. dense entire curves) in the regular part of such Fano fibrations. It immediately implies Theorem \ref{rsurf}.
   
   \begin{proposition}\label{fano-rational}
     Let $S$ be a normal projective surface with only quotient
     singularities, with a non-zero effective $\Q$-divisor $\Delta$ on $S$ such that $K_S+\Delta $  $\Q$-trivial.

     Let $f:S\to B$ be a Fano fibration to a smooth curve $B$
     with relative Picard number $1$. Assume that with $\coeff^h(\Delta) <1$ (resp. $\coeff^h(\Delta)\leq 1)$.

     The following three properties then hold and are equivalent:
      
      1. $\deg(K_B+\Delta_f)<0$ (resp. $\deg(K_B+\Delta_f)\leq 0)$, where $(B,\Delta_f)$ is the orbifold
         basis of $f$.

          2. There exists a rational curve (resp. an entire curve) on $S^{reg}$ which is not $f$-vertical

          3. There exists a $1$-free rational curve through any
         two generic points of $S$ (resp. a dense entire curve)
         which is contained in $S^{reg}$.
 \end{proposition}

   \begin{proof}
     The equivalence of $(1)$, $(2)$ and $(3)$ is due to
     Corollary~\ref{surf'}.
     Property $(1)$ holds because of Lemma~\ref{K_S+D}.
   \end{proof}
   
     \subsection{Examples}
    
    \begin{example}\label{ueno} Let $S:=A/\Z_4$, where $A=E\times E$ and $E:=\C/\Z[\sqrt{-1}]$, and $\Z_4$ is generated by the multiplication by $\sqrt{-1}$ on the two factors simultaneously. Then $S$ (known as ``Ueno surface'')
      has quotient non-canonical singularities and a $\Q$-trivial canonical bundle. If $\pi:S'\to S$ is its minimal resolution, then $-2K_{S'}=\sum E_j$, where $E_j$ are  the $-4$-curves lying over the non-canonical singularities. The preceding Proposition applies to $S'$, but not to $S$ (since all rational curves on $S$ meet some of the non-canonical singularities of $S)$. See \cite{Ca} for a description of the MMP on $S'$. Remark also that $S^{reg}$ contains many dense entire curves (coming from $A)$ although it contains no rational curves.
  \end{example}

 On most normal rational surfaces, every entire curve meets the
  singular locus, as shown (for $m\geq 5)$ in the simplest example below:
  
    \begin{example}\label{rational-ex} Let $f:S\to \P_1$ be a Fano fibration from a normal surface $S$ with Picard number $2$ and with only $2m$ ordinary double points as singularities, these lying in pairs on $m\geq 0$ double fibres of $f$. Such surfaces are obtained by applying the following steps on $m$ distinct fibres, starting from any Hirzebruch surface $S_0$. Blow-up one point, then the intersection point of the two $-1$-curves on the first blow-up. One thus gets a fibre with $3$ components: one double $-1$-curve meeting two reduced $-1$-curves. Then blow-down these last two curves to ordinary double points.  
 
 From Theorem \ref{rsurf} and \cite{CW} (applied as in Lemma \ref{cw}), we get that $S^{reg}$ contains $1$-free rational curves (resp. dense entire curves) if and only if $m\leq 3$ (resp. $m\leq 4)$.       \end{example}

  \section{Dense entire curves transversal to all divisors.}

  One can strengthen the results of \S5 to deal with all divisors $D$
  simultaneously\footnote{This is motivated by Theorem \ref{hp-strong} below.}:

  \begin{theorem}\label{rc-corollary2} Let
    $X$ be a rationally connected complex projective manifold,
    $X$ an analytic subset of codimension at least two
    and $M$ a countable subset of $X\setminus A$.

    Then there exists an entire curve $h:\C\to X\setminus A$
    with $M\subset h(\C)$ such that its image meets any
    irreducible divisor $D\subset X$  transversally
    at some of its smooth points.
  \end{theorem}
    
  \begin{proof}
    Due to theorem~\ref{rc-main} there exists an entire curve
    $h:\C\to X\setminus A$
    with $M\subset h(\C)$. We have to show that this map $h$
    may be chosen in such a way that its image meets any
    irreducible divisor $D\subset X$  transversally
    at some of its smooth points.

    Let $S\subset Chow(X)$ be the subset parametrizing irreducible divisors: it is a countable union of irreducible quasi-projective varieties\footnote{Since $h^1(X,\mathcal O_X)=0$ here, these quasi-projective varieties are Zariski open subsets of projective spaces, and in particular, smooth.}.
      And $S$ is thus the union of an increasing sequence of compact subsets $S_k, k\geq 1$, closures of their interiors, with $S_k$ included in the interior of $S_{k+1}$ for each $k$.
      For each $s\in S$ and $g:\P_1=\C\cup\{\infty\}\to X$, we consider the
      following property $T$
      {\em
        \begin{itemize}
        \item
          There is a complex number $z\in D(1,\frac 12)$ such that
          $c(\P_1)$ intersects $D_s$ transversally in $g(z)$ (i.e.,
          $g'(z)\not\in T_{g(z)}D_s$, in particular
          $g(z)\in D_s^{reg}$).
        \item
          There is no other point $w$ in $\overline{D(1,\frac 12)}$
          with $g(w)\in D_s$.
        \end{itemize}
      }
      This property $T$ is open in $c$ and $s$.
      We may
      a choose a $1$-free rational curve $g_s:\P_1\to X$ such that $g_s(\P_1)$ is transverse to $D_s$ at $q_s=:g_s(1)\in D_s^{reg}$, and has no other zero on the closed disk
      $\overline{D(1, \frac{1}{2})}$.
      In particular, $g_s$ fulfills property $T$ with respect to $s$.
      By the openness of property $T$ there exists an open neighborhood
      $s\in B_s\subset S$
      and a positive number $\alpha_s>0$ such that 
      property $T$ still holds for all $(s',\tilde g)$ where $s'\in B_s$ and
      where
      $\tilde g:\overline{D(1,\frac{1}{2})}\to X$
      is holomorphic such that $d(\tilde g(z), g_{s_m}(z))\leq \alpha_m, \forall z\in \overline{D(1,\frac{1}{2})}$.
      In other words, $\tilde g(\overline{D(1,\frac{1}{2})})$
        meets each $D_{s'}$ ($s'\in B_s$) transversally at a single point
        which is the interior of the disc.

        For each $k>0$, there thus exist finitely many $s\in S_{k+1}\setminus S_k$ such that the open sets $B_s$ cover the (compact) closure $\overline{S_{k+1}\setminus S_k}$ of $S_{k+1}\setminus S_k$.

        We may therefore choose a sequence $s_m$
        (together with a sequence $g_m=g_{s_m}$)  such that
        $S=\cup_{m}B_{s_m}$.

Let $M$ be the set consisting of the given sequence $(x_n)_{n>0}$ of points. We now consider the following sequence of points $y_n, n\geq 1:$ $y_{3m}:=x_m$, $y_{3m+1}:=q'_m$, where $q'_m:=g_m(0)\in g_m(\C)\subset g_m(\P_1)$. Finally, let $y_{3m+2}:=q_m=g_m(1)$, this for every $m>0$. We now construct inductively on $m$ the sequence of rational curves $h_n$ according to this sequence of points $y_n$ as in the proof of Theorem \ref{rc-main}, with the only restriction that for each $m>0$, the curve $h_{3m+2}$ is obtained from $h_{3m+1}$ using the comb defined by $h_{3m+1}$ and $g_m$, that is smoothing this comb. Observe indeed that $h_{m+1}(\P_1)$ contains the point $q'_m\in g_m(\P_1))$, in addition to the points $x_1,\dots, x_m, q'_1,\dots, q'_m, q_1,\dots, q_{m-1}$. 

Choose $h_{3m+2}$ such that $d(g'(z), g_m(z))\leq \frac{\alpha_m}{2^{3m+3}}, \forall z\in \overline{D(1,\frac{1}{2})}$, where $g'(z)=h_{3m+2}(y)$, where $\sigma.y=z$ and $h_{3m+2}:=H({\sigma},.)$ in the notations of the proof of Lemma \ref{rc-lemma}. The transversality of $h_{3m+2}$ to $D_s$ then holds for $s\in B(s,\frac{r_m}{2})$, and also for the subsequent $h_n, n>3m+2$, provided the sequence $\varepsilon_n$ introduced in the proof of Theorem \ref{rc-main} is sufficiently small (it can be inductively choosen). The fact that $h_n$ converges unifomly on compacts of $\C$ to an entire map $h$ satisfying the claims of Theorem \ref{rc-corollary2} is checked as in the proof of Theorem \ref{rc-main}.
    \end{proof}

\subsection{An analogue for $\C^n$.}

Here we discuss dense entire curves in $\C^n$, proving a strong
existence result which will be useful later on.

 \begin{theorem}\label{C^N}

  Let $E,D$ be closed analytic subsets in $\C^n$. Assume
  $\codim(E)\ge 2$, $\codim(D)\ge 1$. Let $S$ be a countable
  subset of $\C^n\setminus E$.

  Then there exists a holomorphic map $F:\C\to\C^n$ such that
  \begin{enumerate}
  \item
    $F(\C)$ is dense in $\C^n$,
  \item
    $S\subset F(\C)$,
  \item
    $F(\C)$ does not intersect $E$, and
  \item
    $F(\C)$ and $D$ are transversal in every point of intersection.
  \end{enumerate}
\end{theorem}

 \begin{proof}
   By enlarging $S$ we may and do assume that $S$ is a dense subset
   of $\C^n\setminus E$.
  We assume that $E$ contains the singular locus of $D$.
  Let $g:\C\to\C^n$ be a holomorphic map with $g(\Z)=S$.
  Let $H$ denote the Fr\'echet vector space of holomorphic maps
  from $\C$ to $\C^n$. For every $v\in H$ we define
  \[
  F_v(z)=g(z)+(e^{2\pi iz}-1)v(z)
  \]
  and observe that $F_v(\Z)=S$.

  For each compact subset $K\subset\C$ we define a subset $W_K$ of
  $H$ by chosing all those $v\in H$ for which the following properties
  hold:
  \begin{enumerate}
  \item
    $F_v(\partial K)\cap D$ is empty,
    \item
      $F_v(K)\cap E$ is empty.
    \item
      $F_v(K)$ intersects  $D$ transversally, i.e.,
      $F'_v(z)\not \in T_{F_v(z)}D$ for $z\in K$ with $F_v(z)\in D$.
  \end{enumerate}

  We claim:
{\em   If $\partial K$ is a smooth real curve,
  then $W_K$ is open and dense in $H$.}

First we explain, how the claim implies the theorem: As a Fr\'echet space,
$H$ has the Baire property. We cover $\C$ by countably many balls $B_i$
and set $W_i=W_{\overline{B_i}}$. Then each $W_i$ is open and dense in
$H$ and due to the Baire property $\cap_i W_i$ is not empty.
Hence we may choose an element $v\in\cap W_i$ and define
$F=F_v$.

Second, we prove the claim.
The topology on $H$ is defined by the sup-Norm. Hence properties
$(i)$ and $(ii)$ are obviously open. Since we discuss holomorphic
functions, locally uniform convergence of functions implies locally
uniform convergence of their derivatives. Therefore it is clear that
$(iii)$ is likewise an open property.

We still have to show density.

We fix a finite-dimensional vector subspace $V$ of $H$ such that
the induced map $V\times\C\to\C^n$ is everywhere submersive.
Then the map
\[
\Phi:V\times (\C\setminus \Z)\to\C^n
\]
given by $(v,z)\mapsto F_v(Z)$ is also everywhere submersive.
It follows that $\Phi^{-1}(D)$ and $\Phi^{-1}(E)$ have at least
codimension $1$ resp.~$2$ in $V\times (\C\setminus \Z)$.
Since $\partial K$ has real dimension one and $K$ has real dimension
$2$, it follows that the projection map onto $V$ maps both
$\Phi^{-1}(D)\cap (V\times\partial K)$ and
$\Phi^{-1}(E)\cap(V\times K)$ to zero measure subsets of $V$.
This implies density for $(i)$ and $(ii)$.

Finally we have to show the density claim for $(iii)$.
We have already seen that the set of all $v\in H$ for which
$F_v(\partial K)\cap D=\emptyset$ holds, is dense. Fix an element $w\in H$
with $F_w(\partial K)\cap D=\emptyset$.  Using the theorem of Rouch\'e, we know
that the degree $deg(F_v^*D|_K$ is constant in a neighbourhood of $w$
in $H$. If $D$ is locally defined by a holomorphic function $h$. then
\[
v\mapsto \Pi_{x\in K\cap F_v^{-1}(D)}dh(F_v'(x)
\]
is holomorphic on $H$, which implies that the complement
of its zero set is
dense.
This completes the proof of the claim and as we have seen before, the
claim implies the theorem.
\end{proof}

\section{Manifolds with $c_1=0$.}\label{c1=0}

We now consider the second main class for testing Conjecture \ref{conjspec} :
compact K\"ahler manifolds $X$  with $c_1(X)=0$. For them, S. Kobayashi already conjectured that their Kobayashi pseudo-distance was identically zero.

{\em Bogomolov decomposition}. A compact K\"ahler manifold with vanishing $c_1$ is (up to finite etale cover)
a direct product of:
\begin{itemize}
\item
  compact complex tori,
\item
  Hyperk\"ahler manifolds, (equivalently: compact K\"ahler manifolds
  which are holomorphically symplectic).
\item
  Calabi-Yau manifolds which are not hyperk\"ahler. They admit a
  nowhere vanishing holomorphic $n$-form ($n$ being the dimension of
  the manifold), but no other non-zero holomorphic differential form.
\end{itemize}
See \cite{B83} and \cite{Bog}.

\subsection{Abelian varieties}

Our result on entire curves in rationally connected manifolds
easily extends to abelian varieties and other compact
complex tori.

\begin{theorem}\label{abelian}
  Let $A$ be a compact complex torus, let $Z$
  be a closed analytic subset of codimension at least two and
  let $M$ be a countable subset of $A\setminus Z$.
    
  Then there exists  an
  entire curve $h:\C \to A\setminus Z$
  with dense image and $M\subset h(\C)$.

  Moreover, every hypersurface $D$
  intersects $h(\C)$ transversally in some point, i.e.,
  there is a point $p=h(t)\in h(\C)\cap D$  with $h'(t)\not\in T_pD$.
  \end{theorem}

\begin{proof}
  We use the universal covering $u:\C^n\to A$.
  Let $E=u^{-1}(Z)$.
  Let $S$ be a dense countable subset
  of $\C^n\setminus E$ with $M\subset u(S)$.
    
  Due to theorem~\ref{C^N} we obtain a holomorphic map
  $F:\C\to\C^n$ with $S\subset F(\C)$ and $F(\C)\cap E=\emptyset$.
  Define $h\stackrel{def}{=}u\circ F$. Now $h:\C\to A$ is a dense
  entire curve which avoids $Z$ and contains $M$ in its image.
  
  Finally observe that due to
  a result from Value Distribution Theory
  (see \cite{NW}, Theorem 6.6.1.)
  for every hypersurface $D$ in $A$ and every Zariski dense entire curve $h$
  there is a point in which $D$ and $h(\C)$ intersect transversally.
 \end{proof}
  
\begin{remark}\label{c_1c_2} Recall that, by \cite{Y}, a compact K\"ahler manifold admits a finite \'etale cover which is a complex torus if and only if it $c_1(X)=c_2(X)=0$.
\end{remark}

\subsection{Manifolds without nontrivial analytic subvarieties}

The following trivial remark nevertheless leads to interesting examples.

\begin{proposition}\label{nosub}
  Let $X$ be a connected normal compact complex space. Assume that $X$ does not contain any non-trivial irreducible complex subvariety (except for $X$ itself and points). If $X$ is not Kobayashi-hyperbolic, then it contains a Zariski-dense entire curve.
\end{proposition}

Indeed, By Brody's Lemma (\cite{Br}, or \cite{Kob}, Theorem 3.6.3),
  there is a non-constant entire curve on $X$,
  which is Zariski-dense, because the whole space $X$ is its only
  positive-dimensional
  analytic subset.

These manifolds are connected to those with $c_1=0$ by means of the following:

\begin{remark}\label{remnosub} A compact K\"ahler manifold $X$ without non-trivial subvariety is conjectured to be either
  a simple compact complex torus
  or to be hyperk\"ahler, in particular, $c_1(X)=0$. See \cite{C06}, Question 1.4, and \cite{CDV}, Conjecture 1.1). This conjecture is proved in dimension $2$ (by Kodaira's classification), and in dimension $3$ in \cite{CDV}, which proves that a smooth compact K\"ahler threefold without subvarieties is a simple torus.\end{remark}

\begin{example} 1.The general\footnote{i.e.: in the countable intersection of Zariski open subsets of the relevant moduli space.} deformation of the Hilbert scheme $Hilb^{m}(K3)$ of $m$ points on a $K3$-surface has no non-trivial
  subvariety, by \cite{V} and is not Kobayashi-hyperbolic, by \cite{V'}.
  Hence Proposition~\ref{nosub} applies and such a manifold does admit a
  Zariski dense entire curve.
  
  2. If $X$ is a compact K\"ahler threefold without subvarieties, it
  is biholomorphic to a compact torus (see \cite{CDV})
  and therefore contains a dense entire curve.
 \end{example}

\begin{remark} Unfortunately our argument gives no information on the
`size' (measured say by its Hausdorff dimension) of the closure of the
entire curve obtained from Brody's Theorem. Hence we can not deduce
property 1.3 or 1.5. or even 1.4. by this method.
\end{remark}
  
\begin{remark}
 The result of \cite{V'} is based on the following construction (\cite{C90}) of entire curves in some Hyperk\"ahler manifolds: if $X$ is a compact K\"ahler Hyperk\"ahler manifold with twistor space $Z$ associated to the Ricci-flat K\"ahler metric on $X$ of a given class $\omega$, then some member $X_s$ of this twistor family contains an entire curve. This entire curve is obtained by deforming the twistor fibres (which are $\P_1's$ with normal bundles direct sums of $\mathcal{O}(1))$ and taking suitable limits. The proof of Verbitsky then shows that this holds in fact for all members of such a family, using the ergodicity of the action of the mapping class group on the Teichm\"uller space and period domain.
\end{remark}

It is quite interesting (but much more difficult) to extend the preceding remarks to the case of `simple' compact K\"ahler manifolds\footnote{These are the main `building blocks' in the construction of arbitrary (non-projective) compact K\"ahler manifolds.}, which are those which are not covered by subvarieties of intermediate dimension, or equivalently, such that their general\footnote{outside a countable union of strict subvarieties.} point is not contained in a strict irreducible compact subvariety. In \cite{C06}, it is conjectured that a `simple' compact K\"ahler manifold is either bimeromorphic to a quotient of torus by a finite group of automorphisms, or is of even complex dimension and carries a holomorphic $2$-form which is symplectic on a nonempty Zariski open subset. In particular: `simple' compact K\"ahler manifolds should have $\kappa=0$, and are of algebraic dimension zero\footnote{that is: have non nonconstant meromorphic function, or equivalently: only finitely many irreducible divisors.}, and are thus `special'. We thus expect them to have dense entire curves. For surfaces, the existence
of Zariski dense entire curves is known by classification: every surface with algebraic dimension zero is simple, and bimeromorphic to either a torus or a $K3$ surface, so far confirming the
conjecture, as we see by the Proposition
  below.

\begin{proposition}\label{a=0,n=2} Let $S$ be a compact K\"ahler surface of algebraic dimension zero. Then $S$ contains a Zariski dense entire curve.
\end{proposition} 

\begin{proof} If $S$ is bimeromorphic to a torus, this is clear. If $S$ is bimeromorphic to a $K3$ surface,we can assume that it is a $K3$ surface. 
In this case, by \cite{BPV},VIII, 3.6, the only connected curves on $S$ are chains of $-2$-curves\footnote{Indeed, all line bundles on $S$ have negative self-intersections, so all nodal classes are classes of $-2$-curves, and by easy computation, these have to meet transversally if not disjoint, never with triple intersections. Although certainly well-known, we do not know a reference.} , and there exists a holomorphic bimeromorphic map $f:S\to S'$ which contracts all curves of $S$ to (singular, normal, Du Val) points of $S'$. Since $d_S\equiv 0$, so is $d_{S'}\equiv 0$, too. Now, $S'$ does not contain nontrivial subvarieties (ie: curves), and so the conclusion follows from Proposition \ref{remnosub}.
\end{proof}

We conclude with the observation that the same holds true in dimension $3$ by the classification of non-projective compact K\"ahler threefolds given in \cite{CHP}.

\begin{proposition}\label{a=0,n=3} Let $X$ be a connected compact K\"ahler threefold of algebraic dimension zero. Then $X$ contains a Zariski dense entire curve.
\end{proposition}

\begin{proof} The claim is obvious, taking Proposition \ref{a=0,n=2} into account, since by \cite{CHP}, \S.9, we get that $X$ is either:

1. `Simple' (hence meromorphically covered by a torus), or bimeromorphically:

2. a $\P_1$-fibration with empty discriminant over a surface of algebraic dimension zero (the surface may be singular, without curves, if $K3$).
\end{proof}

The obstruction to extending the above result to all `special' compact K\"ahler threefolds lies in the cases of `general' projective $K3$ surfaces, normal rational surfaces with quotient singularities and torsion canonical bundle (as the one in Example \ref{ueno}), and Calabi-Yau threefolds.

\subsection{Some Calabi-Yau Manifolds}

These are even more difficult to handle than Hyperk\"ahler manifolds.
In general, we do not know whether the Kobayashi pseudodistance
vanishes. Thus for an arbitrary Calabi-Yau manifold none of the
properties of conjecture \ref{conjspec} is known to be true
(except (3) which states
that $\pi_1$ should be almost abelian).

We just mention $2$ peculiar families for which
we at least know that the Kobayashi pseudo-distance degenerates.
However, also in these cases the other properties (1.2-1.5) are not known.

1. `General' quintics in $\P_4$: they contain $(1,1)$ rational curves of
arbitrarily large degree by \cite{Cl}. This still works in higher dimensions
for general smooth hypersurfaces of degree $(n+2)$ in $\P_{n+1}, n\geq 3$.

2. Double covers of $\P_3$ ramified over a smooth octic: they are covered by elliptic curves by \cite{Vo}, Example 2.17. This still works for double covers of $\P_n$ ramified over a smooth hypersurface of degree $2(n+1)$.

\subsection{Elliptic Fibrations and Elliptic Calabi-Yau Threefolds.}

\begin{proposition} \label{pell} Let $f:X\to B$ be an elliptic  fibration
  from a compact K\"ahler manifold over a rationally connected manifold
  $B$ (i.e.~$f$ is proper, flat and the generic fibers are elliptic curves).
  Assume that $f$ has no multiple fibres in codimension one over $B$.
  Then $X$ contains dense entire curves.
\end{proposition}

\begin{remark}
  With a ``multiple fiber'' we mean a fiber such that {\em every} irreducible
  component of it has multiplicity at least two. 
\end{remark}
  \begin{proof} Let $h:\C\to B$ be any dense entire curve avoiding the subset over which $f$ has multiple fibres. By Theorem \ref{rc-main}, many such entire curves exist. Let
  \[
  Y=E\times_B \C =\{(p,t)\in E\times \C: h(t)=f(p)\}.
  \]
  Now we may regard the natural projection $Y\to\C$ which is
  again an elliptic fibration. Moreover,  there are no multiple fibers
  by construction of $h$.
  Following the arguments of \cite{BL}
  we obtain the existence of a holomorphic map $G:\C^2\to Y$ with
  dense image.
  Because there exists a dense entire curve $\zeta:\C\to\C^2$
  (see e.g.~theorem~\ref{C^N}),
  we are now in a position to define a dense entire curve $g:\C\to E$
  as
  \[
  g(z)\stackrel{def}{=} \pi\left( G(\zeta(z))\right)
  \]
  where $\pi:Y\to E$ is the natural projection map.
\end{proof}

  \begin{theorem}\label{CY3} Let $X$ be a simply-connected Calabi-Yau threefold with terminal singularities and $c_2(X)\neq 0$.
    Assume that $X$ admits an elliptic fibration $f:X\to B$
    over a normal surface $B$.

    Then $X$ contains a dense entire curve.
\end{theorem}

  \begin{proof}
    It follows from \cite{AG}, \cite{O}
    that $B$ is a normal rational surface
 with only quotient singularities and with $-(K_B+D)$ effective for some effective divisor $D$ on $B$ such that $(B,D)$ is Log-terminal, the multiple fibres of $f$
  lying over a finite set $F$.
  From Theorem \ref{rsurf}, we obtain
  the existence of  a dense entire curve inside $(B^{reg}\setminus F)$.
  Now the existence of a dense entire curve in $X$
  follows in the same way as for the preceding Proposition \ref{pell}.
  \end{proof}

\begin{remark} When $c_2(X)=0$, and $X$ is smooth of any dimension, the conclusion still holds since, by \cite{Y}, $X$ is then covered by a complex torus. We thank S. Diverio for this observation.

The conclusion and method of proof of Theorem \ref{CY3} should still apply to Calabi-Yau Manifolds fibered over $\Bbb P_1$ (possibly even an elliptic curve), using the results of \cite{DFM}.\end{remark}

\section{The Nevanlinna version of the Hilbert property}

\subsection{Its statement}

\begin{definition} (\cite{CZ}, \S2.2) Let $X$ be a (smooth) projective variety defined over a number field $k$.
  Then $X$ is said to have the `Weak Hilbert Property' over $k$ (WHP for short)
  \footnote{The classical Hilbert property does not require the covers $Y_j\to X$ to be ramified, it thus implies, by the Chevalley-Weil Theorem, that $X$ to be algebraically simply-connected.}
  if $(X(k)-\cup_{j}\pi_j(Y_j(k)))$ is Zariski-dense in $X$, for any finite set of covers $\pi_j:Y_j\to X$ defined over $k$, each ramified over a non-empty divisor $D_j$ of $X$.
\end{definition}

Note that $X(k)$ being Zariski-dense,
Conjecture \ref{conjspec} implies that
$X$ should  be special, and its fundamental group should be
almost abelian.

In \cite{CZ}, Corvaja-Zannier propose a Nevanlinna version of the $WHP$ in the following form (\cite{CZ}, \S2.4):

\medskip

{\bf Question-Conjecture:} Let $X$ be a simply-connected
smooth projective manifold.
Assume that there exists a Zariski dense entire curve $g:\C\to X$.
For any finite cover $\pi:Y\to X$ ramified over a non-empty divisor,
with $Y$ irreducible, there exists an entire curve $h:\C\to X$ which does
not lift to an entire curve $h':\C\to Y$ (i.e.: such that $\pi\circ h'=h)$.

Let us extend their question to the case of `special' manifolds:

\medskip

{\bf Question-Conjecture:} Let $X$ be a `special' compact K\"ahler manifold. For any finite cover $\pi:Y\to X$ ramified over a non-empty divisor, with $Y$ irreducible, there exists a Zariski dense
entire curve $h:\C\to X$ which does not lift to an entire curve $h':\C\to Y$ (i.e.: such that $\pi\circ h'=h)$.

We shall abbreviate with $NHP(X)$ if $X$ possesses this property, and say that $X$ has $NHP$ (for Nevanlinna-Hilbert Property).

A stronger version is the question whether
there is a single entire curves $h$ on $X$
such that for any finite cover $\pi:Y\to X$ ramified over a
non-empty divisor of $X$, $h$ does not lift to $Y$.
We shall denote with $NHP^+(X)$ this property.

These $NHP$ properties are preserved by finite \'etale covers
and smooth blow-ups.

\begin{lemma}\label{bim-etal} The $NHP$ and $NHP$ for Galois covers properties are preserved by finite \'etale covers, and bimeromorphic equivalence.
\end{lemma}

\begin{proof} Let $f:X'\to X$ be a surjective holomorphic map between connected compact complex manifolds.

Assume first that $f$ is finite \'etale. Let $\pi:Y\to X$ be actually ramified, and $\pi':Y'\to X'$ be deduced from $\pi$ by the base-change $f$. If $h:\C\to X$ is Zariski dense and lifts to $Y$, it lifts to $Y'$ which is \'etale over $Y$. This lifts then lifts some lifting of $h$ to $X'$. A contradiction if $NHP(X')$ holds, and so $NHP(X')$ implies $NHP(X)$. Conversely, let $\pi':Y'\to X'$ be given, actually ramified. If $h':\C\to X'$ is Zariski dense and liftable to $Y'$, then $h:=f\circ h':\C\to X$ is Zariski dense and litable to $Y'$ (which is actually ramified over $X)$, thus $NHP(X)$ is violated as well.

Let us now assume that $f$ is bimeromorphic. Since the fundamental groups of $X$ and $X'$ coincide, so do (up to bimeromorphic equivalence) the covers of $X$ and $X'$ which actually ramify. This (easily) implies the equivalence of $NHP(X)$ and $NHP(X')$, since the existence of lifts of Zariski dense entire curves is also a bimeromorphic property.

The proof for Galois covers is the similar.
\end{proof}

\medskip

A simple tool in finding non-liftable curves is the following:

\begin{proposition}\label{galois-trans}
  Let $h:\C\to X$ be an entire curve and $H$ an hypersurface of $X$ such that there
  exists a regular point $a\in H$ in which $h(\C)$ and $H$ intersect with order of contact $t$.

  Let $\pi:X_1\to X$ be a finite Galois covering with branch locus containing $H$, such that $\pi$ ramifies at order $s\geq 2$ over $H$ at $a$. Then $h$ cannot be lifted to an entire curve $\tilde h:\C\to X_1$ if $t$ does not divide $s$.
  
  In particular, if $h(\C)$ and $H$ meet transversally at $a$, $h$ does not lift to $Y$.
\end{proposition}

\begin{proof}
  Since $\pi$ is Galois, it ramifies at order $s$ at any point of $Y$ over $a\in H$. Since $h(\C)$ intersect at order $s$ at $a$, if it lifted to $Y$, its order of contact with $H$ were a multiple of $s$.\end{proof}

\subsection{ Rationally connected and Abelian manifolds}

We have the following stronger form for
rationally connected manifolds, in which a {\it fixed} entire curve $h$ does not lift to {\it any} Galois ramified cover $\pi:Y\to \P^n$:

\begin{theorem}\label{hp-strong}
  Let $X$ be a rationally connected complex projective manifold
  or a complex compact torus.

Then there exists an entire curve $f:\C\to X$ such that
\begin{enumerate}
\item
The image $f(\C)$ is dense.
\item
$f$ can not be lifted to any ramified Galois covering $\tau:X'\to X$.
\end{enumerate}
\end{theorem}

\begin{proof}
  Combine Theorem~\ref{rc-main}
  resp.~Theorem~\ref{abelian} with Proposition \ref{galois-trans}.
\end{proof}

\subsection{Special surfaces}

\begin{proposition}\label{specsurf}(\cite{Ca04}) Let $S$ be a compact K\"ahler surface. We denote by $S'$ an arbitrary finite \'etale cover of $S$. The following are equivalent:

1. $S$ is special

2. No $S'$ maps meromorphically onto a variety of general type.

3. $\kappa(S)\leq 1$ and $q(S')\leq 2, \forall S'$.

4. $\kappa(S)\leq 1$ and $\pi_1(S)$ is virtually abelian.

5. $S$ is one of the following:

5.1. Rational

5.2. Birational to $\P_1\times E$, with $E$ elliptic.

5.3. Some $S'$ is birational to Enriques, Bielliptic, $K3$, or torus.

5.4. Some $S'$ admits an elliptic fibration $f:S'\to B$ with $B$ either elliptic and no multiple fibre, or $B=\P_1$ and at most $2$ multiple fibres.

\end{proposition}

\begin{theorem}\label{tnhp-surfaces} Let $S$ be a special compact K\"ahler surface. Then $S$ satisfies the Nevanlinna-Hilbert property for Galois covers, except (maybe) if $S$ is a $K3$ surface which is neither elliptic nor Kummer. 

If we assume the Green-Griffiths Lang conjecture, a non-special surface does not fulfill this property.

 \end{theorem}

\begin{proof} We shall first check the property in the cases 5.1-5.3 of Proposition \ref{specsurf}, case 5.4 being treated in Lemma \ref{nhp-surf} below.
If $S$ is rational, it is rationally connected and so satisfies the $NHP$, by Theorem \ref{hp-strong}. If $S$ is in the class 5.2, some birational model has an elliptic fibration, and so the conclusion follows from Lemma \ref{nhp-surf}. Now if $S$ is in the class 5.3, it has a finite \'etale cover bimeromorphic to either an elliptic $K3$ surface, or to a compact torus. In the first case, the conclusion follows frm Lemma \ref{nhp-surf}, in the second from Theorem \ref{hp-strong}.

\end{proof}

\begin{lemma}\label{nhp-surf} A special surface $S$ admitting an elliptic fibration $f:S\to B$ has the $NHP$ for Galois covers.
\end{lemma}
\begin{proof}    

We shall check that $S$ satisfies the condition given by Proposition \ref{galois-trans}. Let $\pi:S'\to S$ be a Galois cover, actually ramified over an irreducible divisor $R\subset S$.
Let $(B,D_f)$, with $D_f:=\sum_j(1-\frac{1}{m_j}).\{b_j\}$ be the orbifold base of the fibration $f$, where $b_j\in B,\forall j$, and $m_j$ is the multiplicity\footnote{The `classical' and `non-classical' versions coincide for elliptic fibrations.} of the fibre $S_{b_j}$ of $f$ over $b_j$.

Since $S$ is special, $deg(K_B+D_f)\leq 0$. For any $b\in B$, there thus exists an entire map $h:\C\to (B,D_f)$ , which goes through $b$, and which is an orbifold entire curve (ie: such that its order of contact with any $b'\in B$ is equal to the multiplicity of $D_f$ at $b'$, which is equal to $1$ if $b'$ is not any one of the $b_j's$).

Let $f_h:S_h:=S\times_B\C\to \C$ be deduced form $f$ by the base change $h:\C\to B$. Then $f_h$ has no multiple fibre, and (by \cite{BL}) admits a holomorphic section $s:\C\to S_h$ going through any given reduced component $F$ of any of the fibres of $f_h$.

The main result of \cite{BL} then implies the following:

\begin{lemma} Assume that $R$ is either a reduced component of some fibre of $f_h$, or that $f(R)=B$. There exists then a holomorphic map $H:\C^2\to X$ which is unramified over its image, which is a Zariski open subset of $S$ meeting $R$. 

Composing $H$ with an suitable injection $j:\C\to \C^2$ with dense image and meeting transversally $H^{-1}(R)$ at some point, Proposition \ref{galois-trans} and Lemma \ref{bim-etal} imply $NHP(S)$ for Galois covers.
\end{lemma}

We still need to deal with the case when $R$ is a non-reduced component of some fibre $F:=f_h^{-1}(0)$ of $f_h$. Let then $m>1$ be the multiplicity of $R$ in $F$. Let $\mu:\C\to \C$ be defined by $\mu(z)=z^m$, let $k:=h\circ\mu:\C\to B$, and  let $f_k:S'_k\to \C$ be deduced from $f_h$ by the base-change $\mu:\C\to C$ after taking a smooth bimeromorphic model $S'_k$ of $S_k:=S\times_B\C$.

The natural map generically finite map $\sigma_k:S'_k\to S_h$ is not \'etale, but it is \'etale over the generic point of $R$, and also outside the fibre $F$. Applying again \cite{BL} to the fibration $f_k:S'_k\to \C$, we obtain a holomorphic map $H_k:\C^2\to S'_k$ which is unramified over the generic point of the inverse image $R_k$ of $R$, and thus a Zariski dense holomorphic map $h_k:\C\to S'_k$ which is transversal to $R_k$. The map $\sigma_k\circ h_k:\C \to S_h$ is thus transversal to $R$ at some of its generic points, establishing $NHP(S)$ for Galois covers. \end{proof}

\subsection{Removing the Galois condition}

\begin{theorem}
  Let $X$ be a complex projective manifold. Assume that $X$ is rationally
  connected or that there exists a surjective and submersive
  holomorphic map
  $\rho:\C^N\to X$.

  Let $\pi:Y\to X$ be a finite map with non-empty ramification.

  Then there exists a dense entire curve $h:\C\to X$
  which can not be lifted to $Y$ (i.e., there is no entire curve
  $\tilde h:\C\to Y$ with $h=\pi\circ\tilde h$.)
\end{theorem}

This result differs from theorem~\ref{hp-strong} in two points:
First, the ramified covering $\pi$ is no longer required to be
Galois. Second (this is the price we pay for dropping the Galois
condition), here the non-liftable dense entire curve $h$ may depend
on $\pi$.

\begin{proof}
  First we discuss the case $\dim_{\C}X=1$. Then $X$ is a rational or elliptic curve.
  If $X$ is an elliptic curve, $Y$ is of genus at least  $2$ and therefore hyperbolic, i.e.,
  there is no non-constant holomorphic map from $\C$ to $Y$.
  If $X$ is a rational
  curve, we chose the natural injection $\C\subset\P_1$ as $h:\C\to\P_1$.
  Using the theorem of Liouville one may show easily: For every non-constant
  holomorphic map $\tilde h$ from $\C$ to a compact Riemann surface $Y$
  the image must be dense. This implies that such a map $\tilde h$
  can not be a lift of $h$.

  Now we start the proof of the general case, i.e., $\dim_{\C}X\ge 2$.
  Let $R\subset Y$ denote the ramification divisor and
  let $B=\pi(R)\subset X$ be the branching locus.
  Fix $p\in X\setminus B$. For every $q\in\pi^{-1}(p)$
  we choose a real continuous curve
  $c_q:[0,1]\to Y$ with $c_q(0)=q$, $c_q(1)\in R$ and
  $\pi(c_q(t))\not\in B$ for $0\le t<1$. Since $\dim_{\C}X\ge 2$,
  we may, by slightly perturbing the curves $c_q$ if necessary,
  assume that for $q\ne q'$ we have
  $\pi(c_q([0,1]))\cap\pi(c_q'([0,1]))=\{p\}$.
  Thus $T=\cup_q \pi(c_q([0,1]))$ is a tree,
  i.e., a simply-connected real one-dimensional simplicial complex
  which may be visualized as $d=\#\pi^{-1}(p)$ line segments glued
  together in one point.
  Due to Proposition~\ref{treeCn} resp.~Proposition
  \ref{treeRC}  we obtain an entire curve
  $h_n:\C\to X$ with embeddings $\zeta_n:T\to\C$ such that
  $\lim_{n\to\infty}h_n\circ\zeta_n=i$ where $i:T\to X$
  is the inclusion map.
  For sufficiently large $n$ we now may replace $T$ by $(h_n\circ\zeta_n)(T)$
  and obtain an obstruction to the lifting of the entire curve $h_n$
  by using proposition~\ref{tree-obstruction}.

  Finally, it is a consequence of Corollary~\ref{cor-dense-ini} that
  the entire curve may be chosen in such a way that it has dense image.
\end{proof}

\begin{proposition}\label{tree-obstruction}
  Let $\pi:Y\to X$ be a finite ramified covering of complex  manifolds
  with ramification divisor
  $R\subset Y$ and branching locus $B=\pi(R)$.

  Let $p\in X\setminus B$. For every $q\in \pi^{-1}(p)$ let $c_q:[0,1]\to Y$
  be a continuous real curve with  $c_q(0)=q$ and $c_q(1)\in R$
  and $\pi(c_q(t))\not\in B$ for all $0\le t <1$.

  Let $h:\C\to X$ be an entire curve
  such that for every $q\in\pi^{-1}(p)$ there is a smooth real curve
  $\gamma_q:[0,1]\to\C$ with $\gamma_q(0)=0$,
  $h\circ\gamma_q=\pi\circ c_q$
  and $h'(\gamma_q(1))\not\in TB$.

  Then there does not exist a lift $\tilde h:\C\to Y$, i.e., there
  is no holomorphic map $\tilde h:\C\to Y$ with $\pi\circ\tilde h=h$.
\end{proposition}

\begin{proof}
  Assume the converse. Since $\gamma_q(0)=0$ for all $q$, we have
  \[
  h(0)=h(\gamma_q(0))=\pi(c_q(0))= \pi(q)=p.
  \]
  Fix $q\in\pi^{-1}(p)$ such that $\tilde h(0)=q$.
  Observe that
  \begin{itemize}
  \item
    $\pi\circ\tilde h\circ\gamma_q= h\circ\gamma_q = \pi\circ c_q$,
  \item
    $\tilde h\circ\gamma_q(0)=q=c_q(0)$,
  \item
    $\pi\circ c_q(t)\not\in B$ for all $t<1$.
  \end{itemize}
  It follows that $\tilde h\circ\gamma_q= c_q$ and therefore
  $(\tilde h\circ\gamma_q)(1)\in R$.
  In combination with $h'(\gamma_q(1))\not\in TB$ this yields
  a contradiction.
\end{proof}
    
\begin{proposition}\label{treeCn}
Let $T$ be a (real) tree, i.e., a simply-connected finite one-dimensional
  simplicial complex.

  Let $X$ be a complex projective manifold with a surjective
  and submersive
  holomorphic map $\rho:\C^N\to X$.

  Let $c:T\to X$ be a continuous map.
  Let $\epsilon>0$ and let $d(\ ,\ )$ be a distance function
  on $X$ induced by a hermitian metric.
  Then there exists an entire
  curve $h:\C\to X$ and a continuous map 
  $\tilde c:T\to\C$ such that
  $h(\C)$ is dense in $X$ and
  $d(c(x),h(\tilde c(x))<\epsilon$
  for all $x\in T$.
\end{proposition}
\begin{proof}
  We chose a closed embedding $\tilde c:T\to\C$ and an infinite
  discrete subset $S\subset\C$ with $S\cap\tilde c(T)=\{\}$.
  Upon replacing $c$ by a small deformation, we may assume that
  $c(T)$ contains no critical values of $\rho$.
  Since $\rho$ is surjective, we may lift $c:T\to X$ to a continuous map
  $\hat c:T\to\C^N$. The theorem of Arakelyan/Nersesyan (\cite{Ne})
  implies that
  continuous maps from $\tilde c(T)\cup S$ to $\C^N$ may uniformly
  be approximated by holomorphic maps from $\C$ to $\C^N$.
  Hence
  the continuous map $\hat c\circ(\tilde c)^{-1}:\tilde c(T)\to\C^N$
  may be approximated uniformly by holomorphic maps $\tilde h$
  from $\C$ to $\C^N$ in such a way that $\tilde h(\C)$ is dense in $\C^N$.
  This implies the assertion, by taking $h\stackrel{def}=\rho\circ\tilde h$.
\end{proof}

\begin{proposition}\label{treeRC}
  Let $T$ be a (real) tree, i.e., a simply-connected finite one-dimensional
  simplicial complex.

  Let $X$ be a rationally connected complex projective manifold.

  Let $c:T\to X$ be a continuous map.
  Let $\epsilon>0$ and let $d(\ ,\ )$ be a distance function
  on $X$ induced by a hermitian metric.
  Then there exists a free
  rational curve $h:\P_1\to X$ and a continuous map 
  $\tilde c:T\to\P_1$ such that $d(c(x),h(\tilde c(x))<\epsilon$
  for all $x\in T$.
\end{proposition}

\begin{proof}
  We decompose the tree into  arcs. Locally, each arc admits
  an approximation by free rational curves, see lemma~\ref{lem-arc}.
  Due to compactness, we obtain a decomposition into finitely any arcs
  such that each arc admits an approximation by a free rational curve.
  Using comb smoothing (in the form given by corollary~\ref{cor-comb-tree})
  we obtain an approximation
  by one free rational curve.
\end{proof}
  
\begin{lemma}\label{lem-arc}
  Let $X$ be a rationally connected projective manifold with a
  distance function $d(\ ,\ )$ (induced by some hermitian metric).
  Let $\epsilon>0$, $r\in\N$ and let $\gamma:[-1,1]$ be a regular real
  curve (i.e.~$\gamma$ is $C^1$ with $\gamma'(t)\ne 0\ \forall t$).

  Then there exists a positive number $0<\delta<1$ such that for
  every $-\delta<a<b<\delta$ there is an $r$-free rational curve
  $c:\P_1\to X$ satisfying
  \begin{enumerate}
  \item
    $d(c(t),\gamma(t))<\epsilon$ for every $t\in [a,b]$.
  \item
    $c(a)=\gamma(a)$, $c(b)=\gamma(b)$,
  \item
    $c'(a)=\gamma'(a)$, $c'(b)=\gamma'(b)$.
  \end{enumerate}

\end{lemma}

\begin{proof}
  Without loss of generality we assume $r\ge 2$.
  We start by choosing an $r$-free rational curve $h:\P_1\to X$
  with $h(0)=\gamma(0)$ and $h'(0)=\gamma'(0)$.
  We consider the deformations of $h$.
  Since $h$ is free, these deformations are unobstructed. Thus we may
  fix a ball $B\subset V=H^0(\P_1,h^*TX)$
    parametrizing deformations
  of $h$
  by a map $\Phi:B\times\P_1\to X$.
  
  Choosing the ball small enough, we may assume
  $d(\Phi(p,t),h(t))<\epsilon/2$ for all $p\in B$, $t\in[-1,+1]$.
  We may furthermore assume that $\Phi$ extends to $\bar B\times\P_1$
  where $\bar B$ is the closure in $V$ and compact.

  Because $h$ is 1-free, for every pair $(a,b)$, $-1\le a<b\le+1$ 
  the evaluation map
  \[
  B\ni p \mapsto \Phi(p,a)\times \Phi(p,b)
  \]
  contains an open neighborhood of $(a,b)$ in $X\times X$ in its image.
  Using compactness, we choose $\eta>0$ such that for all $a,b\in[-1,1]$
  and all $a',b'\in X$ with $d(\gamma(a),a')<\eta$, $d(\gamma(b),b')<\eta$
  there exists a parameter $p\in B$ with $\Phi(p,a)=a'$ and $\Phi(p,b)=b'$.

  Next we chose $\delta>0$ in such a way that
  \[
  d(c(t),\gamma(t))< \min\{\eta,\epsilon/2\}\ \forall t, |t|<\delta.
  \]

  Given any two numbers $a,b$ satisfying $-\delta<a<b<\delta$,
  we have $d(c(a),\gamma(a))<\eta$ and $d(c(b),\gamma(b))<\eta$ and
  therefore
  may choose $p\in B$ with $\Phi(p,a)=\gamma(a)$
  and $\Phi(p,b)=\gamma(b)$. Then we define $c(t):=\Phi(p,t)$.
  We obtain  a $r$-free rational curve $c$. Evidently $c$ satisfies
  property $(ii)$ of the assertion.
  Furthermore, $c$ satisfies $(i)$ because
  \[
  d(c(t),\gamma(t))\le d(c(t),h(t))+d(h(t),\gamma(t))
  \le \epsilon/2 + \epsilon/2=\epsilon.
  \]

  Finally we note that $(iii)$ may be achieved in the same way
  as in
  Corollary~\ref{rc-corollary1}.
\end{proof}

  \begin{remark}
  As pointed out to us by F.~Forstneric, there is a related result
  by A.~Gournay, who proved- under a certain additional `regularity'
  hypothesis- a Runge type approximation theorem (see \cite{Go}) 
  on domains in Riemann surfaces mapped to compact almost complex manifolds containing `free' rational curves.
  Where applicable, the statement of the above lemma
  can be deduced from Gournay result as follows:
  First approximate a given map $\gamma:[0,1]\to X$
  by a real-analytic map, then extend it to a holomorphic map
  $g$ on some simply-connected
  open neighborhood of $[0,1]$. If a certain `regularity' condition is
  fulfilled,  Gournays theorem then implies that
  $g$ may be approximated by a rational curve, implying the statement of
  the above lemma. The `regularity' assumption, although satisfied for `generic' almost complex structures,
   seems however difficult to check on a given complex structure.
\end{remark}

\begin{lemma}
  Let $T$ be a tree, and let $c:T\to C$ be a continuous closed embedding
  with $C=\{(z,w)\in\C^2:zw=0\}$.
    
    Then there exists
     a sequence of non-zero complex numbers
    $\epsilon_n$
    with $\lim_{n\to\infty}\epsilon_n=0$ and 
  a sequence of continuous maps
    $c_n:T\to C_{\epsilon_n}$ such that $\lim c_n=c$ 
    uniformly on $T$.
\end{lemma}
\begin{proof}
  The assertion is easily verified if $c(T)$ is contained in one
  of the two irreducible components of $C$. Hence we assume that
  $c(T)$ is not contained in one of the irreducible components,
  implying that $c(T)$ contains $(0,0)$, since this is the only
  point in which the two components intersect.
  
  Let $0\in T$ denote the point with $c(0)=(0,0)$.
  Let $T_i^{\circ}$ denote the connected components of $T\setminus\{(0)\}$
    and let $T_i=T_i^{\circ}\cup\{(0)\}$ be its closure.
    
  We choose an injective continuous curve $\gamma:[0,1]\to\C$ with
  $\gamma(0)=0$ such that $\gamma(t)+x^2\ne 0$ for all $t\in]0,1]$
      and $(0,x)$, $(x,0)\in T$.

      Since $T$ is a tree, we know that $T_i\times[0,1]$ is simply-connected.
      The choice of $\gamma$
  allows us thus to choose a branch  of $\sqrt{\gamma(t)+x^2}$
  (with $(x,0)$ resp.~$(0,x)$ in $T_i$ and $t\in[0,1]$). We choose the
  branch $h_i(x,t)$  with $\sqrt{0+x^2}=x$ for some point
  $(x,0)$ resp.~$(0,x)$ in $T_i$. By continuity it follows that
  $\sqrt{0+x^2}=x$ for all $x$. Consequently $\lim_{t\to 0}
  h_i(t,x)=x$ for all $(x,0)$ resp.~$(0,x)$ in $T_i$.
  We define $c_{t,i}:T_i\to\C^2$ as
  \[
  c_{t,i}(x)= \frac 12\left( h_i(t,x)+x, h_i(t,x)-x\right)
  \]
  if $T_i\subset\C\times\{(0)\}$ and
  \[
  c_{t,i}(x)= \frac 12\left( h_i(t,x)-x, h_i(t,x)+x\right)
  \]
  if $T_i\subset\{(0)\}\times\C$.
  We observe that
  \[
  \lim_{t\to 0} c_{t,i}=c.
  \]
  Now $h_i(t,0)^2=\gamma(t)$. We fix a branch of the square
  root along $\gamma$:
  \[
  \xi(t)\stackrel{def}=\sqrt{\gamma(t)},\ t\in[0,1]
  \]

  Then for each $i$ there is an element $\sigma_i\in\{+1,1\}$
  such that
  \[
  h_{t,i}(0)=\sqrt{\gamma(t)+0}=\sigma_i\xi(t)
  \]
  and consequently
  \[
  c_{t,i}(0)=(\sigma_i\xi(t),\sigma_i\xi(t))
  \]
  The problem is that for distinct $i,j$ we may have $\sigma_i\ne\sigma_j$.

  We enumerate the $T_i$ as $T_1,\ldots, T_r$. Since $T$ is a tree, and
  hence contains no closed loops, each subtree $T_i$ contains a unique
  edge ending at $(0,0)$.

  Now we define $c_t:T\to\C^2$ recursively as follows:
  If $c_t$ is defined on $\cup_{i\le k}T_i$, then we choose a small
  curve inside $C_{\gamma^2(t)}$
  starting at $c_t(0)$ and ending at $c_{t,k+1}(0)$. Because
  \[
  ||(\xi(t),\xi(t)-(-\xi(t),-\xi(t))||=2\sqrt{ 2\gamma(t)}
  \]
  and $\lim_{t\to 0}\gamma(t)=0$, this small curve can be chosen
  as small as desired if $t$ is small enough.
  Thus we may attach this curve to $c_{t,k+1}(T_{k+1})$ without
  changing this subtree to much.

  In this way we can construct a family of continuous maps
  $c_t:T\to C_{\gamma^2(t)}$ with $\lim_{t\to 0}c_t=c$
  (uniformly on $T$).
\end{proof}

\begin{corollary}\label{cor-comb-tree}
  Let $C$ be a comb in a rationally connected projective manifold,
  let $T$ be a tree
  and let $c:T\to C$ be a closed continuous  embbeding.

  Then $c$ may be approximated by $h_n\circ c_n$ where $h_n$ are free
  rational curves and $c_n:T\to\P_1$ are closed continuous embeddings.
\end{corollary}

\section{Dense entire curves of order zero}

For an entire curve $f:\C\to X$ with values in a K\"ahler manifold $X$
the characteristic function (in its Ahlfors-Shimizu form) is defined
as
\[
T_f(r)=\int_0^r \left( \int_{\Delta_t}f^*\omega \right) \frac{dt}{t}
\]
where $\omega$ is a K\"ahler form on $X$ and $\Delta_t=\{z\in\C:|z|<t\}$.

The {\em order} $\rho_f$ of $f$ is defined as
\[
\rho_f=\limsup_{r\to\infty} \frac{\log T_f(r)}{\log r}
\]

We will show that the entire curves in rationally connected manifolds
constructed above in this article may be constructed in such a way that
the order $\rho_f$ equals zero.

This is interesting in view of the following facts:

\begin{itemize}
\item
  For every rational map $f:\C\to X$ we have $T_f(r)=O(\log r)$ and
  consequently $\rho_f=0$.
\item
  For an abelian variety $A$ and a non-constant entire curve $f:\C\to A$
  we have $\rho_f\ge 2$.
\item
  The notion of the order $\rho_f$ may generalized to non-degenerate
  holomorphic maps from $\C^n$ to $X$. If (for $n=\dim(X)$) there exists
  a non-degenerate holomorphic map $f:\C^n\to X$ with $\rho_f<2$,
  then $X$ is rationally connected (see \cite{CW"},\cite{NW2}).
\end{itemize}

We observe that from the definition of the characteristic function
given above the following may be easily deduced:

{\em If $f_n$ is a sequence of holomorphic maps from $\C$ to some
K\"ahler manifold $(X,\omega)$ which converges locally uniformly
on $\C$ to a holomorphic map $f:\C\to X$, then the sequence
of characteristic functions $T_{f_n}:[1,\infty[\to\R^+$ converges
    locally uniformly on $[1,\infty[$ to $T_f:[1,\infty[\to\R^+$.}
          
Let $X$ be a projective rationally connected manifold. We fix a K\"ahler class
$\omega$. For any rational curve $f:\P_1\to X$ the degree
$\deg(f)$ is defined
as the volume of its image with respect to the K\"ahler form, i.e.,
\[
\deg(f)=\int_{\P_1}f^*\omega
\]
For any $c>0$ let $F_c$ denote the set of all free rational curves of
degree at most $c$.
Due to the K\"ahler assumption the degree is invariant under deformations.
Hence there is a constant $U>0$ such that for any two given points $p,q$ there
is a free rational curve $f$ in the family $F_U$ connecting $p$ with $q$.

Let $g:[1,\infty[\to [1,\infty[$ be a continuous function with the following
    properties
    \begin{enumerate}
    \item
      $r\mapsto \frac{g(r)}{\log r}$ is monotone increasing
      on $[2,\infty[$
      and unbounded.
    \item
      $\lim_{r\to\infty}\frac{\log g(r)}{\log r}=0$.
    \item
      $g(2)\ge 2U \log 2$
    \end{enumerate}

    For example, we may take
    $g(r)=\min\{ (\log r)^2, 2 U\log r\}$, or: $g(r)=\min\{ (\log(\log r)).\log r, 2 U\log r\}$.
    
    We choose recursively a sequence $f_n$ of rational curves in X,
    with $\deg(f_n)=nU$.

    We start with a rationcal curve $f_1$ of degree $\le U$.
    For $r\ge 1$ we have:
           \begin{align*}
        T_{f_{1}}(r) &\le  \int_1^{r}
        \left(\int_{D_t}f_{1}^*\omega\right) \frac{dt}{t} \\
        &\le U\log r\\
       \end{align*}
           Since $g(2)\ge 2 U\log 2$, we have $T_{f_1}(2)\le \frac 12 g(2)$.
           By the monotonicity of  $g(r)/\log r$ it follows
           that
           \[
           T_{f_1}(r)\le \frac 12 g(r)\quad \forall r\ge 2
           \]

           We claim that the recursive construction of the rational
           curves $f_n$ as made above may be carried out in such a way
           that
 \[
    T_{f_n}(r)\le \left( 1 - \frac{1}{2^{n+1}}\right) g(r)\quad\forall r\ge 2, \forall n,
      \]
           
    Given $f_n$
    with \[
    T_{f_n}(r)\le \left( 1 - \frac{1}{2^{n+1}}\right) g(r)\quad\forall r\ge 2,
      \]
     first we choose $R_n>2$ such that
    \[
    \frac{g(R_n)}{\log R_n}> 2(n+1)U
    \]

    Now we choose $f_{n+1}$ such that on $\{z:|z|\le R_n\}$
    it is close enough to $f_n$ in order to ensure
    \[
    T_{f_{n+1}}(r)\le \left( 1 - \frac{1}{2^{n+2}}\right) g(r)
      \quad\forall r\le R_n
      \]

      Let us now consider $T_f(r)$ for $r>R_n$:
      \begin{align*}
        T_{f_{n+1}}(r) &= T_{f_{n+1}}(R_n) + \int_{R_n}^r
        \left(\int_{D_t}f_{n+1}^*\omega\right) \frac{dt}{t} \\
        &\le T_{f_{n+1}}(R_n) + \int_{R_n}^r
         (n+1)U \frac{dt}{t} \\
&=         T_{f_{n+1}}(R_n) +
        (n+1)U \left( \log r -\log R_n\right) \\
      \end{align*}
        By the choice of $R_n$, the condition $r\ge R_n$ and the
        monotonicity of $g(r)/\log r$ we have $2U(n+1) < g(r)/\log r$.
        Therefore
        \[
        T_{f_{n+1}}(r) \le T_{f_{n+1}}(R_n) + \frac 12 g(r)- (n+1)U\log R_n
        \]
        Because of $\deg(f_{n+1})=(n+1)U$ we have
       \begin{align*}
        T_{f_{n+1}}(R_n) &\le  \int_1^{R_n}
        \left(\int_{D_t}f_{n+1}^*\omega\right) \frac{dt}{t} \\
        &\le (n+1)U\log R_n\\
       \end{align*}
       Hence
       \[
       T_{f_{n+1}}(R_n) - (n+1)U\log R_n<0
       \]
       and consequently
       \[
        T_{f_{n+1}}(r) \le \frac 12 g(r)\quad  \forall r\ge R_n
        \]
        It follows that
       \[
       T_{f_{n+1}}(r) \le \left( 1 - \frac{1}{2^{n+2}}\right)
       g(r)\quad  \forall r\ge 2
        \]

        Finally, we obtain an entire curve $f:\C\to X$ as
        $f=\lim f_n$ and our construction implies
        \[
        T_f(r)\le g(r)
        \]
        which implies that $f$ is of order zero:
        \[
        \rho_f=\limsup_{r\to\infty}\frac{\log T_f(r)}{\log r}
        \le \limsup_{r\to\infty}\frac{\log g(r)}{\log r}=0
        \]
        The entire curves so constructed have thus growth order zero.

\end{document}